\documentclass[12pt, a4paper]{asl}

\title[Large cardinals and morasses]{Large cardinals and gap-1 morasses}

\author{Andrew D. Brooke-Taylor}
\revauthor{Brooke-Taylor, Andrew D.}

\address{Department of Mathematics\\
University of Bristol\\
University Walk\\
Bristol, BS8 1TW, UK}

\email{Andrew.Brooke-Taylor@bristol.ac.uk}

\thanks{Andrew Brooke-Taylor was supported for this research by Austrian
Science Fund (FWF) project P16790-N04.}

\author{Sy-David Friedman}
\revauthor{Friedman, Sy-David}

\address{Kurt G\"odel Research Center for Mathematical Logic\\
University of Vienna\\
W\"ahringer Stra\ss e 25\\
A-1090, Wien, Austria}

\email{sdf@logic.univie.ac.at}

\thanks{Sy-David Friedman was supported by FWF projects P16790-N04 and
\mbox{P19375-N18.}}

\newcommand{\al}{\alpha}
\newcommand{\be}{\beta}
\newcommand{\ga}{\gamma}
\newcommand{\de}{\delta}
\newcommand{\ze}{\zeta}
\newcommand{\io}{\iota}
\newcommand{\ka}{\kappa}
\newcommand{\la}{\lambda}
\renewcommand{\phi}{\varphi}
\newcommand{\si}{\sigma}

\newcommand{\Power}{\mathcal{P}}
\newcommand{\st}{\,|\,}
\newcommand{\compat}{\,\|\,}
\newcommand{\forces}{\Vdash}
\newcommand{\sat}{\vDash}

\renewcommand{\P}{\mathbb{P}}
\newcommand{\U}{\mathbb{U}}
\newcommand{\Iff}{\longleftrightarrow}
\renewcommand{\iff}{\leftrightarrow}
\renewcommand{\implies}{\rightarrow}

\newcommand{\into}{\hookrightarrow}
\newcommand{\restricted}{\upharpoonright}
\newcommand{\restr}{\!\upharpoonright\!}

\newcommand{\Exists}{\,\exists}

\newcommand{\cf}{\textrm{cf}}

\DeclareMathOperator{\dom}{dom}

\DeclareMathOperator{\Fn}{Fn}
\newcommand{\GCH}{\textrm{GCH}}

\newcommand{\id}{\textrm{id}}

\newcommand{\Lim}{\textrm{Lim}}
\newcommand{\ND}{\textrm{ND}}
\newcommand{\ORD}{\textrm{Ord}}
\newcommand{\ot}{\textrm{ot}}

\newcommand{\rlm}{\textrm{rlm}}

\DeclareMathOperator{\supp}{supp}
\newcommand{\trcl}{\textrm{trcl}}

\newcommand{\val}{\textrm{val}}
\newcommand{\ZF}{\textrm{ZF}}
\newcommand{\ZFC}{\textrm{ZFC}}

\newcommand{\calD}{\mathcal{D}}

\newcommand{\calL}{\mathcal{L}}
\newcommand{\LSTG}{\calL_{STG}}

\newcommand{\calS}{\mathcal{S}}

\newcommand{\mcequiv}[1]{\stackrel{#1}{\sim}}
\newcommand{\morl}{\shortmid\!\prec}
\newcommand{\closermorl}{\shortmid\negmedspace\prec}
\newcommand{\notmorl}{\shortmid\!\nprec}
\newcommand{\tal}{\theta_\alpha}
\newcommand{\tla}{\theta_\lambda}
\newcommand{\pixy}{\pi_{xy}}

\newcommand{\HkaP}{H_{\ka^+}}
\newcommand{\HlaP}{H_{\la^+}}

\newtheorem{thm}{Theorem}
\newtheorem{defn}[thm]{Definition}
\newtheorem{lemma}[thm]{Lemma}

\newtheorem{prop}[thm]{Proposition}

\renewcommand{\qedsymbol}{$\dashv$}

\newcommand{\emin}{\emph}

\begin{document}

\begin{abstract}
We present a new partial order for directly forcing morasses to exist 
that enjoys a
significant homogeneity property.  We then use this forcing in a reverse 
Easton iteration to obtain an extension universe with morasses at every
regular uncountable cardinal, while preserving all $n$-superstrong
($1\leq n\leq\omega$), hyperstrong and 1-extendible cardinals.  In the latter
case, a preliminary forcing to make the GCH hold is required.
Our forcing yields morasses that satisfy an extra property related to 
the homogeneity of the partial order; we refer to them as 
\emph{mangroves} and prove that their existence is equivalent to the 
existence of morasses.  Finally, we exhibit a partial order that forces
universal morasses to exist at every regular uncountable cardinal, and use
this to show that universal morasses are consistent with $n$-superstrong,
hyperstrong, and 1-extendible cardinals.
This all contributes to the second author's \emph{outer model programme},
the aim of which 
is to show that $L$-like principles can hold in outer models which 
nevertheless contain large cardinals.
\end{abstract}

\maketitle

\section{Introduction}

Morasses are combinatorial structures formulated by Jensen to abstract out 
properties 
of $L$ useful for proving cardinal transfer theorems.
Originally thought to be
such complex objects that one should not attempt
to understand them outside the context of their construction in $L$
(according to \cite{Sta:VRe}),
they were opened up to broader study and deeper understanding by
Velleman \cite{Vell:MDF}, who amongst other things presented a 
particularly elegant partial order with which one can force morasses to exist.
This lends itself to use in the second author's
\emph{outer model programme} in which, as a counterpoint to the 
well-known inner model programme, the goal is to force to obtain 
$L$-like outer models for various large cardinals.
Indeed, in \cite{SDF:LCL}
it is shown that one may force morasses to exist at every regular 
uncountable cardinal, while preserving a given $n$-superstrong or 
hyperstrong cardinal.

In this article we extend these results.  In addition to $n$-superstrong and
hyperstrong cardinals we consider 1-extendible cardinals, for which
a careful consideration of forcing the GCH as a preliminary step is appropriate.
Further, we exhibit a new forcing partial order to obtain the morasses,
which allows us to preserve \emph{all} cardinals of the considered kinds.
This is achieved by modifying Velleman's forcing to give it a degree of
homogeneity.
Finally, we show how one may force universal morasses to exist.
This new forcing again fits in with the techniques of \cite{SDF:LCL},
and so we obtain that universal morasses can be forced while preserving 
an $n$-superstrong, hyperstrong or 1-extendible cardinal.  
This answers a question of Pereira to the authors.

The reader may be aware that the existence of morasses 
has been shown by Velleman \cite{Vell:SiM} to 
be equivalent to the existence of much simpler objects called
\emin{simplified morasses}. 
However, while these latter structures are much more
manageable, 
the only known way to obtain them in $L$ is via morasses themselves.
Thus, to stay more directly
in keeping with the goal of forcing $L$-like models, we stick with 
morasses.  On the other hand, we do use Velleman's \cite{Vell:MDF}
axiomatisation of morasses, which already strips away some of the detail
of morasses as defined by Jensen.  This will have ramifications when 
we force to obtain universal morasses --- see section \ref{univMorasses} below.

\section{Preliminaries: the GCH and 1-Extendible Cardinals}\label{GCH}

Our forcing constructions will rely on having the GCH hold in the ground model,
so we discuss here how to attain it by a preliminary forcing.
The iterated forcing we use to obtain the GCH is a fairly standard,
natural one, 
at each stage collapsing the old
$2^{\aleph_\al}$ to $\aleph_{\al+1}$.  
The main issue is to show that forcing with this
particular partial order does preserve the large cardinals we are interested 
in.  In many cases, this was achieved in \cite{SDF:LCL}; 
we show here that 
similar arguments go through for 1-extendible cardinals.

Following \cite{Kun:ST}, we denote by $\Fn(I,J,\la)$ the set of 
partial functions from $I$ to $J$ with cardinality less than $\la$.
\begin{defn}
The \emin{GCH Partial Order}\index{partial order!GCH} $P$ is the
reverse Easton iteration of $\langle\dot Q_\al\st\al\in\ORD\rangle$,
in which direct limits are taken only at strongly inaccessible cardinals and 
inverse limits are taken at other limit stages,
where $\dot Q_\al$ is the 
canonical $P_\al$-name for 
$\Fn(\aleph_{\al+1},2,\aleph_{\al+1})$, with $\aleph_{\al+1}$ to
be evaluated in $V[G_\al]$.
\end{defn}

For the rest of this section $P=P_{\ORD}$ will denote the GCH Partial Order,
with $P_\al$ denoting the iteration after $\al$ stages, 
$G=G_{\ORD}$ denoting a $P$-generic over $V$,
$G_\al$ denoting $G\restr P_\al$,
and $Q_\al$ denoting 
$\Fn(\aleph_{\al+1},2,\aleph_{\al+1})$
in $V[G_\al]$.

Before proceeding with the analysis of this iteration,
let us recall the basic properties of the
forcing poset $\Fn(\aleph_{\al+1},2,\aleph_{\al+1})$.
It is $\aleph_{\al+1}$-closed
and $(2^{\aleph_\al})^+$-cc,
so the only cardinals collapsed by forcing with 
it 
are those $\ka$ such that
$\aleph_{\al+1}<\ka\leq 2^{\aleph_\al}$.
Of course, since $\Fn(\aleph_{\al+1},2,\aleph_{\al+1})$
is equivalent to $\Fn(\aleph_{\al+1},2^{\aleph_\al},\aleph_{\al+1})$,
a surjection $\aleph_{\al+1}\to 2^{\aleph_\al}$ is added, 
and all such $\ka$ are indeed collapsed.
By a nice names argument, the continuum function is unchanged by this 
forcing at and above the ground model $2^{\aleph_\al}$, 
and $2^{\aleph_{\al+1}}$ in the extension will be the ground model 
$2^{2^{\aleph_\al}}$.
By $\aleph_{\al+1}$-closure, the continuum function is unchanged below
$\aleph_\al$, and after the forcing we have 
$2^{\aleph_\al}=\aleph_{\al+1}$.

Now, to the GCH partial order itself.  To demonstrate that ZFC is preserved
by this class forcing, we show that it is tame (see \cite{SDF:FSCF}).
\begin{lemma}\label{GCHForcTame}
The GCH Partial Order is tame.
\end{lemma}
\begin{proof}
By the Factor Lemma (see for example \cite{Jech:ST}, Lemma 21.8), 
$P$ may be written as
$P_\al*\dot P^{\al}$, where $\dot P^{\al}$ names
the iteration starting with $Q_\al$.  Each iterand in $\dot P^{\al}$
is forced to be $\aleph_{\al+1}$-closed, and inverse limits are taken
everywhere except at inaccessibles of $V$, which will remain inaccessible
after forcing with $P_\al$ (indeed, it can be shown by induction that
$P_\al$ has a dense suborder of size at most $\beth_{\al+1}^V$). 

We therefore have that for all $\al$
$$
1_{P_\al}\forces\dot P^\al\text{ is }\aleph_{\al+1}\text{-closed},
$$
and this implies that $P$ is tame --- see \cite{SDF:FSCF}, Lemmas 2.22 and
2.31.
\end{proof}

\begin{lemma}\label{GCHForcGCH}
Let $\al$ be an ordinal or $\ORD$.
Then for every ordinal $\ga\leq\al$, 
$\beth_\ga^{V[G_\al]}=\aleph_\ga^{V[G_\al]}$.
\end{lemma}
\begin{proof}
As in the proof of Lemma \ref{GCHForcTame} above, $P_\al$ may be factorised as
$P_\ga*\dot P^{[\ga,\al)}$, where $\dot P^{[\ga,\al)}$ names
the iteration starting with $Q_\ga$, and we have
$$
1_{P_\ga}\forces\dot P^{[\ga,\al)}\text{ is }\aleph_{\check\ga+1}\text{-closed}.
$$
It hence suffices to show that 
$\beth_\ga^{V[G_\ga]}=\aleph_\ga^{V[G_\ga]}$
for every $\ga$.  This may be proven by induction 
on $\ga$, with the successor step following from the 
discussion of $\Fn(\aleph_{\al+1},2,\aleph_{\al+1})$ above.
\end{proof}

In particular, these two lemmas combine to show that forcing with $P$
yields a model of $\ZFC+\GCH$.


In \cite{SDF:LCL}, it is shown that when forcing with the GCH forcing
$P$, any given hyperstrong or $n$-superstrong cardinal may be preserved. 
In fact, the same proof essentially gives the following stronger 
statement.

\begin{thm}
Forcing with the GCH partial order $P$ preserves all 
hyperstrong cardinals and $n$-superstrong cardinals, for all $n\in\omega+1$.
\end{thm}
\begin{proof}
Observe that each iterand in the GCH forcing is very homogeneous:
given a condition $p\in Q_\al=\Fn(\aleph_{\al+1},2,\aleph_{\al+1})$ 
and a $Q_\al$-generic $g$, 
there is an 
automorphism $\pi$ of $Q_\al$ such that
$p\in\pi``g$.  
Such strong homogeneity is known to hold for an entire iteration if
it is forced to hold for each iterand; 
see \cite{DoF:HIM1C} for the details.
It follows that all the arguments for Theorem~4.2 of \cite{SDF:LCL},
which in particular choose $P$-generics containing a certain master
condition for each large cardinal $\ka$, 
can be carried out within $V[G]$.
Therefore, all
large cardinals of the listed kinds are preserved in $V[G]$.
\end{proof}

We shall show that the same result is true for 1-extendible cardinals.
Recall their definition:

\begin{defn}
A cardinal $\ka$ is \emph{1-extendible} if there is a $\la>\ka$ and
an elementary embedding $j:H_{\ka^+}\to H_{\la^+}$ with critical point $\ka$.
\end{defn}
Equivalently, $\ka$ is 1-extendible if and only if there is a $\la>\ka$ and
an elementary embedding $V_{\ka+1}\to V_{\la+1}$ with critical point $\ka$.
The above definition will be more convenient for our purposes, however, 
as $H_{\ka^+}$
forms a model of $ZFC$ minus the Power Set axiom, 
which we shall refer to as 
$ZFC^-$.  1-extendible cardinals lie in consistency strength between 
superstrong and hyperstrong cardinals (see \cite{Kan:THI}, Proposition 26.11)
and play an interesting role in the outer model programme, as
Cummings and Schimmerling \cite{CuS:ISq} have shown them to be essentially the
strongest large cardinals compatible with the principle $\square$.

\begin{thm}\label{GCH1Ext}
For any model $V$ of $\ZFC$, 
there is a class-generic extension $V[G]$ of $V$ such that $V[G]\sat\GCH$,
and in which every 1-extendible cardinal of $V$
remains 1-extendible.
\end{thm}
\begin{proof}
We shall of course force with the GCH forcing $P$.
Let $\ka$ be 1-extendible in $V$, let $j:H_{\ka^+}\to H_{\la^+}$ in $V$
witness this fact, and let $G$ be $P$-generic over $V$.  
We shall show that $j$ 
may be lifted to an embedding 
$j^*:H_{\ka^+}^{V[G]}\to H_{\la^+}^{V[G]}$
witnessing the 1-extendibility of $\ka$ in $V[G]$.

The GCH forcing $P$ may be factorised as $P_\ka*P^\ka$,
with $P^\ka$ a $\ka^+$-closed forcing.  Hence,
$H_{\ka^+}^{V[G]}=H_{\ka^+}^{V[G_\ka]}$.
Similarly, 
$H_{\la^+}^{V[G]}=H_{\la^+}^{V[G_\la]}$.
Now because $\ka$ is inaccessible, $P_\ka$ has cardinality $\ka$ and
lies in $H_{\ka^+}$.
Recall that forcing over a model of $\ZFC^-$ yields a model of $\ZFC^-$,
and note that
$G_\ka$ is generic over $H_{\ka^+}$ for the set forcing $P_\ka$.
We claim that in fact 
$H_{\ka^+}^{V[G_\ka]}=H_{\ka^+}^V[G_\ka]$, the generic extension of the 
model $H_{\ka^+}^V$ of $\ZFC^-$ by $G_\ka$.
For this we need to show that if $\sigma$ is a $P_\ka$-name in $V$ such that
$\sigma_{G_\ka}\in H_{\ka^+}^{V[G_\ka]}$, then there is a name 
$\tau\in H_{\ka^+}^{P_\ka}$ such that $\sigma_{G_\ka}=\tau_{G_\ka}$.
But now every element of $H_{\ka^+}$ can be obtained as the
Mostowski collapse of a relation on $\ka$ given by a subset of $\ka\times\ka$
(and this process does not appeal to the Power Set Axiom).
Every subset of $\ka\times\ka$ in the extension $V[G_\ka]$
has a nice name of the form 
$$
\bigcup_{(\al,\be)\in\ka\times\ka}\{\check{(\al,\be)}\}\times A_{(\al,\be)}
$$
where each $A_{(\al,\be)}$ is an antichain in $P_\ka$.
Since $|P_\ka|=\ka$, such nice names lie in $H_{\ka^+}^V$.
Therefore, every subset of $\ka\times\ka$ in $V[G_\ka]$ is 
also in $H_{\ka^+}[G_\ka]$, and consequently we indeed have that
$H_{\ka^+}^{V[G_\ka]}=H_{\ka^+}^V[G_\ka]$.
Of course by elementarity, $\la$ is also inaccessible, and so also
$H_{\la^+}^{V[G_\la]}=H_{\la^+}^V[G_\la]$.

We are thus reduced to showing that $j:H_{\ka^+}\to H_{\la^+}$
may be lifted to an elementary embedding 
$j^*:H_{\ka^+}[G_\ka]\to H_{\la^+}[G_\la]$.
We wish to apply 
Silver's technique of lifting embeddings
(see \cite{Cum:IFE}, Section 9),
and must show that 
$j``G_\ka\subseteq G_\la$.  But now $j\restr P_\ka$ is the identity function,
so $j``G_\ka=G_\ka\subset G_\la$.  Therefore 
$j$ indeed lifts to $j^*:H_{\ka^+}^{V[G]}\to H_{\la^+}^{V[G]}$,
and so $\ka$ is 1-extendible in $V[G]$.
\end{proof}

In Section \ref{Mangrove1Extend}, we will show that 1-extendible cardinals
$\ka$ may be preserved while forcing morasses to exist.  To show that the
forcing is pretame from the perspective of $H_{\ka^+}$,
we will first force GCH as above, 
and use the fact that after our GCH forcing, $V[G]=L[G]$, 
and further $H_{\ka^+}^{V[G]}=L_{\ka^+}[A]$, where $A$ is 
taken to be a class predicate over $H_{\ka^+}$ for the Cohen set
added at stage $\ka$. 
This gives a stratification of our model
$H_{\ka^+}$ of $\ZFC^-$ into sets, 
as is required for the proof that the forcing relation 
is definable for a pretame class forcing --- see \cite{SDF:FSCF}, 
Theorem 2.18.
Only then will we be able to conclude that
the forcing relation for our forcing is
definable.

However, $p\forces\varphi(\sigma_0,\ldots,\sigma_{n-1})$ 
will only be definable relative to $A$, and so
the usual lifting lemma will not suffice to lift the embedding $j$
witnessing 1-extendibility, as $j$ need not respect arbitrary class predicates.
We show here that in fact $j$ \emph{does} respect the predicate $A$
given by the generic, giving a mild strengthening of Theorem 
\ref{GCH1Ext}.
The crux of the proof will be the definability of the forcing relation
for $P_{\ka+1}$, which from the point of view of $H_{\ka^+}$ is a
class forcing.  
The usual argument for pretame class forcing (Theorem 2.18 of \cite{SDF:FSCF}) 
is not applicable,
since we don't yet have a stratification of $H_{\ka^+}$ into sets.
Nevertheless, we will be able to 
demonstrate the definability of forcing for $P_{\ka+1}$ 
because it is such a well-behaved forcing partial order.

\begin{defn}\label{LSTG}
The language $\calL_{STG}$\index{$\calL_{STG}$}
is the language obtained from the language of set theory
by adding a single unary predicate $G$.
\end{defn}

\begin{lemma}\label{GCHForcingHkaPDefinable}
Let $\ka\in M$ be an inaccessible cardinal,
and let $Q_\ka$ denote $\Fn(\ka^+,2,\ka^+)$, the forcing at stage
$\ka$ in the GCH forcing above.
Then the forcing relation $\forces\phi$ 
for forcing over $\HkaP$ with $Q_\ka$
is definable over $\HkaP$ for $\calL_{STG}$ formulae $\phi$,
where the predicate $G$ is to be interpreted as the $Q_\ka$-generic,
and the Truth Lemma holds: a formula is true in the extension if and only if
it is forced by some condition in the generic.
\end{lemma}
\begin{proof}
We follow the expository style of \cite{Kun:ST}, defining a relation 
$\forces^*\phi$ for each $\phi$, which will then be seen to be equivalent to
$\forces\phi$.
We will make use of the closure properties of $Q_\ka$ to make this
definition for atomic formulae, 
from which point it may be extended to all formulae as usual.
As a first step, we give a definition for forcing the predicate $G$ to hold:
define $p\forces^*G(\sigma)$ to mean
\begin{equation}
\forall q\leq p\Exists\,r\leq q\Exists s\geq r(r\forces^*\sigma=\check s),
\label{ClassForcGDefn'}
\end{equation}
that is, it is dense below $p$ to force $\sigma$ to be something specific 
greater than yourself.
If $p\in G$, then by genericity there is some $r\in G$ and some 
$s$ greater than or equal to $r$ and consequently also in $G$ such that
$r\forces^*\sigma=\check s$, whence $\sigma_G\in G$.  Conversely, suppose
$\sigma_G\in G$ for every $G\ni p$;  let $q\leq p$ and fix some $G\ni q$. 
Once we have shown that the forcing relation is definable for
atomic formulae, it will follow 
from the truth lemma there is some $r'\in G$ such that
$$r'\forces\sigma=\check{(\sigma_G)}.$$  Taking $r\in G$ such that
$r\leq q$, $r\leq r'$, and $r\leq\sigma_G$, we have that
$r$ makes (\ref{ClassForcGDefn'}) hold with $\sigma_G$ as $s$.
Hence, (\ref{ClassForcGDefn'}) is indeed a formal definition encapsulating
the statement that for every generic containing $p$, $\sigma_G$ is in $G$,
and it is clear that the argument extends to prove the Truth Lemma valid
for formulae involving the predicate $G$.

With forcing the predicate $G$ so defined, we can move on to show that
the forcing relation is indeed definable for atomic formulae.  
Let $\sigma$ and $\tau$
be $Q_\ka$-names in $\HkaP$, and let  
$$
R_{\sigma,\tau}=
\trcl(\{\sigma,\tau\})\cap P,
$$ 
the set of conditions hereditarily appearing in either
$\si$ or $\tau$.
Note that since the names $\si$ and $\tau$ are in $\HkaP$, 
$|R_{\sigma,\tau}|\leq\ka$.
Let $d_{\si,\tau}=\bigcup_{r\in R_{\si,\tau}}\dom(r)$; 
then $d_{\si,\tau}$ is a subset of $\ka^+$ of
size at most $\ka$.

Now, for conditions $q$ such that $d_{\si,\si}\subseteq\dom(q)$,
we may recursively define the evaluation of $\si$ at $q$ by
$$
\val(\si,q)=\{\val(\rho,q)\st
\exists X\subseteq\dom(q)(\langle\rho,q\restr X\rangle\in\sigma)\}.
$$
Observe that if $G$ is any $Q_\ka$-generic over $\HkaP$ containing such a $q$, 
then $\si_G=\val(\si,q)$.  
Likewise, for any $G$ which is $Q_\ka$-generic over $\HkaP$ and 
any element $q$ of $G$ with domain containing $d_{\si,\si}$, we have
$\val(\si,q)=\si_G$.  Thus, the evaluation of any name in the generic 
extension may be entirely determined by a single condition.

With this in hand, we may define
$$
p\forces^*\si\in\tau\quad\Iff\quad
\forall q\leq p\big(\dom(q)\supseteq d_{\si,\tau}\implies
\val(\si,q)\in\val(\tau,q)\big)\phantom{.}
$$
and
$$
p\forces^*\si=\tau\quad\Iff\quad
\forall q\leq p\big(\dom(q)\supseteq d_{\si,\tau}\implies
\val(\si,q)=\val(\tau,q)\big).
$$
Then clearly $p\forces^*\sigma\in\tau$ if and only if $\si_G\in\tau_G$ for 
every $Q_\ka$-generic $G\ni p$, 
and $p\forces^*\si=\tau$ if and only if $\si_G=\tau_G$ for
every $Q_\ka$-generic $G\ni p$.
Moreover, if $\si_G\in\tau_G$, then any $p\in G$ with 
$\dom(p)\supseteq d_{\si,\tau}$ will force $\si=\tau$, and the Truth Lemma 
follows.
\end{proof}

We are now ready for our strengthening of Theorem \ref{GCH1Ext}.
Let use denote by $G(\ka)$ the generic from stage $\ka$ of $G$, that is,
the $Q_\ka$-generic associated to $G$.

\begin{thm}\label{GCH1ExtG}
Let $V$ be a model of $\ZFC$ and let $\ka$ be a 1-extendible cardinal in
$V$ with 1-extendibility witnessed by an elementary embedding 
$j:H_{\ka^+}\to H_{\la^+}$.  Let $V[G]$ be a $P$-generic extension of $V$.
Then there is a $G'\subset V[G]$ which is $P$-generic over $V$ such that
$V[G]=V[G']$ and the lift $j^*$ 
of $j$ to $H_{\ka^+}^{V[G']}$ 
(as in the proof of Theorem \ref{GCH1Ext}) is elementary between the
$\calL_{STG}$-structures
$\langle H_{\ka^+}^{V[G']},G'(\ka)\rangle$ and
$\langle H_{\la^+}^{V[G']},G'(\la)\rangle$. 
\end{thm}
\begin{proof}
In Theorem~\ref{GCH1Ext} the lift $j^*$ was constructed in
$\HkaP^{V[G_\ka]}$, after observing that $P^\ka$ is $\ka^+$-closed and
hence does not affect $\HkaP$.  
We now claim that the
$P^\ka$ generic $G^\ka$ may be chosen so that this same $j^*$ is
also elementary for formulae in the language $\LSTG$.
By the Truth Lemma, every $\LSTG$ sentence $\phi$ true in 
$\langle\HkaP^{V[G_{\ka+1}]},G(\ka)\rangle=
\langle\HkaP^{V[G]},G(\ka)\rangle$ is forced (over $\HkaP^{V[G_\ka]}$)
to be true by some $p\in Q_\ka$.  
Now by the definability of the forcing relation, we have
$$
p\forces_{Q_\ka}\phi(\si_0,\ldots,\si_n)\iff
j^*(p)\forces_{Q_\la}\phi(j^*(\si_0),\ldots,j^*(\si_n)).
$$
Therefore, if $j^*``G(\ka)\subseteq G(\la)$, we will have that $j^*$ is
elementary from 
$\langle H_{\ka^+}^{V[G]},G(\ka)\rangle$ to
$\langle H_{\la^+}^{V[G]},G(\la)\rangle$. 
But now $|G(\ka)|^{V[G_\la]}=(\ka^+)^{V[G_\la]}<\la$,
so $j^*``G(\ka)$ is a pairwise-compatible 
set of conditions in $\Fn(\la,2,\la)^{V[G_\la]}$
(indeed, in $\Fn(\la,2,\la)^{V[G_\ka]}$) of size less than $\la$.
Moreover, this set lies in $V[G_\la]$, since $j$, $G_\ka$ and $G(\ka)$ all do.
Hence, $m=\bigcup j^*``G(\ka)$ is a single condition in
$\Fn(\la,2,\la)^{V[G_\la]}=Q_\la$, and so if $G$ is chosen such that this 
master condition\index{master condition} lies in $G(\la)$, then
we will indeed have $j^*``G(\ka)\subset G(\la)$.
This is easy to arrange by modifying $G(\la)$ to obtain $G'(\la)$.
Let $q$ be the element of $G(\la)$ with $\dom(q)=\dom(m)$, and define
$\psi:Q_\la\to Q_\la$ changing partial functions
by switching their values between $0$ and $1$ on 
points where $q$ and $m$ disagree, that is,
$$
\psi(f)(\al)=\begin{cases}
0&\text{if }q(\al)\neq m(\al)\text{ and }f(\al)=1\\
1&\text{if }q(\al)\neq m(\al)\text{ and }f(\al)=0\\
f(\al)&\text{otherwise}
\end{cases}
$$
for all $f\in Q_\la$ and $\al\in\dom(f)$.
Clearly $\psi$ is an (involutive) automorphism of $Q_\la$ in $V[G]$, whence
$G'(\la)=\psi``G(\la)$ is $Q_\la$-generic over $V[G_\la]$,
and $V[G_{\la+1}]=V[G_\la*G'(\la)]$.  
The ``tail'' generic $G^{\la+1}$ will of course still be
$P^{\la+1}$-generic over $V[G_\la*G'(\la)]$.  Hence,
considering $P$ as $P_\la*Q_\la*P^{\la+1}$, we may take
$G'=G_\la*G'(\la)*G^{\la+1}$, which is $P$-generic over $V$ and
satisfies $V[G]=V[G']$.
Moreover, $G'(\la)$ contains the master condition $m$ determined by
$G(\ka)=G'(\ka)$, so $j^*``G'(\ka)\subset G'(\la)$, and so
$j^*$ is elementary from
$\langle H_{\ka^+}^{V[G']},G'(\ka)\rangle$ to
$\langle H_{\la^+}^{V[G']},G'(\la)\rangle$, as required. 
\end{proof}

\section{Basic definitions}\label{morasses}

The definition for morasses that we use will be that of 
Velleman \cite{Vell:MDF}.
In particular, we retain the notation M.1--M.7 for axioms of a morass, used 
both there and in Devlin \cite{Dev:Con}.
We also follow Velleman in separating out the following subsidiary
definition, although we change the terminology slightly to emphasise the 
order preservation property.

\begin{defn}
A function $\pi:\al\to\beta$ between two ordinals is a
successor, limit, zero and order preserving (SLOOP) function if
\begin{itemize}
\item for all $\gamma<\al$, $\pi(\gamma+1)=\pi(\gamma)+1$, and
\item for all limit 
$\gamma<\al$, $\pi(\gamma)$ is also a limit, and
\item $\pi(0)=0$, and
\item for all $\gamma<\delta<\al$, $\pi(\gamma)<\pi(\delta)$.
\end{itemize}
\end{defn}
Note in particular that the composition of SLOOP functions will yield a 
SLOOP function.

Morasses will be defined based on a set $\calS$ of ordered pairs of ordinals.
For such an ordered pair $x=\langle\al,\beta\rangle$, let us denote 
$\al$ by $l(x)$ (the \emin{level} of $x$), and $\beta$ by $o(x)$
(the \emin{order} of $x$). 
Also, note that we here use the word ``tree'' in the liberal sense where others
might use ``forest'': our tree will have many root nodes.

\begin{defn}\label{morass}
For any uncountable regular cardinal $\ka$, a \emph{$(\ka, 1)$-morass}
(or simply \emin{morass} when $\ka$ is clear from the context)
consists of:
\renewcommand\theenumi {\roman{enumi}}
\begin{enumerate}
\item a subset $\calS$\index{$\calS$!of a morass} of 
$(\ka\times\ka)\cup(\{\ka\}\times\ka^+)$, and
\item A tree order $\morl$ on $\calS$, and
\item For every pair $\langle x,y\rangle$ 
of elements of $\calS$ with $x\morl y$, 
a function $\pi_{xy}:o(x)+1\to o(y)+1$,
\end{enumerate}
\renewcommand\theenumi {\arabic{enumi}}
such that the following conditions hold.  
\begin{description}
\item[Left-alignment] For each $\al\leq\kappa$, let 
$\tal=\{\beta\,|\,\langle\al,\beta\rangle\in\calS\}$.
Then $\tal$ is in fact an ordinal.  Moreover, 
$\theta_\ka=\ka^+$, and for $\al<\ka$, $0<\tal<\ka$.
\item[Monotonicity] For $x$ and $y$ in $\calS$, $x\morl y$ implies $l(x)<l(y)$.
\item[Commutativity] If $x\morl y\morl z$ then 
$\pi_{yz}\circ\pi_{xy}=\pi_{xz}$.

\item[M.1] For each pair $x,y\in\calS$ with $x\morl y$, 
$\pi_{xy}$ is an SLOOP function, and $\pi_{xy}(o(x))=o(y)$.

\item[M.2] Suppose $x\morl y\in\calS$ and $\nu<o(x)$.  Let 
$w=\langle l(x), \nu\rangle$ and $z=\langle l(y), \pi_{xy}(\nu)\rangle$.  
Then
$w\morl z$ and $\pi_{wz}=\pi_{xy}\restricted(\nu+1)$.

\item[M.3] For all $y\in\calS$, $\{l(x)\st x\morl y\}$ is closed in $l(y)$.

\item[M.4] For all $y\in\calS$, if $o(y)+1\neq\theta_{l(y)}$ 
(as defined for \emin{Left-alignment} above), then
$\{l(x)\st x\morl y\}$ is unbounded in $l(y)$.  In particular, if $\al$ is a 
successor ordinal, then $\tal=1$.

\item[M.5] For all $y\in\calS$, if 
$\{l(x)\st x\morl y\}$ is unbounded in $l(y)$, then
$o(y)=\bigcup\{\pi_{xy}\!``o(x)\st x\morl y\}$.

\item[M.6] Suppose $x\morl y\in\calS$ and $o(x)$ is a limit ordinal.  
Let $\nu=\sup(\pi_{xy}\!``o(x))$ and let $z=\langle l(y),\nu\rangle$.
Then $x\morl z$ and $\pi_{xz}\restricted o(x)=\pixy\restricted o(x)$.

\item[M.7] Suppose $x\morl y\in\calS$, $l(x)<\al<l(y)$, 
$o(x)$ is a limit ordinal, and
$o(y)=\sup(\pixy\!``o(x))$.  If
$$\forall\nu<o(x)\,\exists\gamma\,(\langle\al,\gamma\rangle\morl
\langle l(y),\pixy(\nu)\rangle),$$ 
then there is a $\gamma$ such that
$\langle\al,\gamma\rangle\morl y$.

\end{description}
\end{defn}

Note in particular axiom M.5, in the case where $l(y)=\kappa$: for any 
$\tau<\ka^+$, we have $\tau$ expressed as the (increasing, by commutativity) 
union of the images of the maps 
$\pi_{x\langle\ka,\tau\rangle}\restricted o(x)$ for 
$x\morl\langle\ka,\tau\rangle$.  
Further note that this isn't just a variant of 
$\ka$ many things ``adding up'' to $\ka$, and then $\tau$ being bijective with 
$\ka$, but something more direct: the maps $\pi_{x\langle\ka,\tau\rangle}$ are
order preserving, and so the ordinals mapping into $\tau$ must to some extent 
reflect the structure of $\tau$.  
For example, if $\tau$ is a successor ordinal and
so has a largest element, then 
for $x\morl\langle\ka,\tau\rangle$ with $l(x)$ sufficiently large, $o(x)$ must
also have a largest element.
Jensen's orignal definition in fact called for even stronger preservation
properties; see \cite{Dev:Con} for details.
Also observe that for any $x\morl y$ in $\calS$, $o(x)\leq o(y)$
since by M.1, $\pixy$ is a strictly order-preserving function from
$o(x)+1$ to $o(y)+1$.

When we force morasses to exist in the presence of multiple large cardinals, 
it will be convenient for the sake of preservation of the large cardinal 
property to use a partial order which lends itself to homogeneity arguments 
(see Section \ref{homogeneity} below).  
The upshot will be that the generic morass satisfies
a useful property not possessed by the morasses obtained by forcing with 
Velleman's partial order defined in \cite{Vell:MDF}.
We give here a name for the kind of morass that we obtain.

\begin{defn}\label{mangal}
Suppose $M=\langle\calS,\morl,\langle\pixy\rangle_{x\morl y}\rangle$ is a 
$(\ka,1)$-morass.  An ordinal $\al<\ka$ is a \emin{mangal} of $M$ if for all 
$x, y\in\calS$ with $x\morl y$ and $l(x)<\al<l(y)$, there is a 
$z\in\calS$ with
$l(z)=\al$ such that $x\morl z\morl y$.  If the set of mangals of $M$ is cofinal
in $\ka$, we say that $M$ is a \emph{$\ka$-mangrove}, or simply 
\emin{mangrove} when $\ka$ is clear from the context.
\end{defn}

An immediate question is whether the existence of a morass implies the
existence of a mangrove.  It turns out that the answer is yes, as we shall
show in Section \ref{ManfrMor}, using Velleman's theorem that the
existence of a morass is equivalent to a certain forcing axiom.

Another natural question is what the set of mangals of a morass
can look like.  A first result in this direction is the following.

\begin{prop}\label{MangalsClosed}
For any morass $M$, the set of mangals of $M$ is closed in $\ka$. 
\end{prop}
\begin{proof}
Suppose that $\ga<\ka$ is a limit point of the set of mangals of a
morass $M=\langle\calS,\morl,\pi\rangle$, and let $x, y\in\calS$ satisfy
$x\morl y$ and $l(x)<\ga<l(y)$.  For each mangal $\al$ of $M$ such that
$l(x)<\al<\ga$, there is a $z_\al$ with $l(z_\al)=\al$ and
$x\morl z_\al\morl y$.  But then by M.3, there must be some $z_\ga$ with
$l(z_\ga)=\ga$ and $x\morl z_\ga\morl y$; thus, $\ga$ is also a mangal of $M$.
\end{proof}
In particular, if $M$ is a mangrove, the set of mangals of $M$
is a closed unbounded subset of $\ka$.

\section{Forcing morasses to exist}\label{MangroveForcing}

As is the case for many combinatorial structures, 
one can force with a partial order of partial morasses to get a morass in the 
generic extension, as we shall show below.
We use a partial order similar to those described in \cite{SDF:LCL} and
\cite{Vell:MDF}.  However, ours will differ in that we strengthen the 
requirements for an extension of a condition, ensuring that each
condition of the generic goes up to a mangal of the ultimate morass.
Naturally, this will make our generic morass a mangrove, whereas it is not hard
to check that a generic morass for Velleman's partial order in \cite{Vell:MDF}
will be far from a mangrove; indeed,
it will be $\ka$-branching at the node $\langle 0,0\rangle$, and 
so cannot have any mangals at all!
However, 
many of the details of our proof will remain essentially the same as in 
that paper.  In particular, we repeatedly use the basic
construction given there for extending morass conditions, tweaking
it to fit the particular requirements in each case.

\subsection{Definitions}

Let $x\morl_i y$ denote that $y$ is an immediate $\morl$-successor of $x$.

\begin{defn} \label{MorassCondition}
For any uncountable regular cardinal $\ka$, a 
\emph{$(\ka, 1)$-morass condition}\index{morass condition}
consists of:
\renewcommand\theenumi {\roman{enumi}}
\begin{enumerate}
\item a subset $\calS$\index{$\calS$!of a morass condition} 
of $((\lambda+1)\times\ka)\cup(\{\ka\}\times\ka^+)$ for
some $\lambda<\ka$, and
\item A tree order $\morl$ on $\calS$, and
\item For every pair $\langle x,y\rangle$ 
of elements of $\calS$ with $x\morl y$, 
a function $\pi_{xy}:o(x)+1\to o(y)+1$,
\end{enumerate}
\renewcommand\theenumi {\arabic{enumi}}
such that 
\begin{enumerate}
\item \label{MorCondLAS}
\emph{Left-alignment} holds for all $\al\leq\lambda$, and the set
$S=$\hbox{$\{\beta\,|\,\langle\ka,\beta\rangle\in\calS\}$} 
\index{$S^p$ (morass condition level $\ka$)}
contains $0$ and is closed
under ordinal successors and predecessors.
\item \label{MorCondBij}
Let $f$ be the order-preserving bijection from $\ot(S)$ to $S$.
Then $\ot(S)\leq\theta_\lambda$, and for each $\nu<\ot(S)$,
$\langle\lambda,\nu\rangle\morl\langle\ka,f(\nu)\rangle$.
\item \emph{Monotonicity}, \emph{Commutativity}, \emph{M.1},
\emph{M.2}, \emph{M.3} and \emph{M.5} hold.  
Axiom \emph{M.4}
holds for those $y\in\calS$ such that $l(y)\leq\lambda$.
Axioms \emph{M.6} and \emph{M.7} 
hold for those $x$ and $y$ such that $x\morl_iy$ and 
$l(y)\leq\lambda$.
\end{enumerate}
\end{defn}

For our analysis, we will need to extend the new notion of a mangal
to also be applicable to morass conditions.
\begin{defn}\label{MorCondMangal}
Let $p=\langle\calS,\morl,\pi\rangle$ be a morass condition.
An ordinal $\al\leq\la$ is a \emin{mangal} of $p$ if for all 
$x, y\in\calS$ with $x\morl y$ and $l(x)<\al<l(y)$, there is a 
$z\in\calS$ with
$l(z)=\al$ such that $x\morl z\morl y$.  
\end{defn}

Our definition of a morass condition
is the same as that used in \cite{Vell:MDF}; 
our partial order will differ in the definition of $\leq$.  
For this reason, we continue to refer to these conditions as 
\emph{morass} conditions, even though we shall call our partial order the
\emph{mangrove} forcing.
Also note that, since it arises frequently and the meaning is fairly clear, 
we shall consistently abuse notation, writing $x\in p$ to mean that
$x$ is an element of the set $\calS$ for the morass condition $p$.

A few comments about the definition of morass conditions are in order 
at this point.
Note that requiring Axiom M.2 when $l(y)=\ka$ implicitly entails imposing 
the condition on $S$ and the various $\pixy$ that 
for any $x\morl y$ with $l(y)=\ka$,
$\pixy``(o(x)+1)\subset S$.
Requiring M.7 for $x\morl_iy$ is the same as positing the
non-existence of $\al$ fitting the antecedent of that axiom, since
$x\morl\langle\al,\ga\rangle\morl y$ contradicts $x\morl_iy$.

Observe that by requirement \ref{MorCondLAS} above, 
$\ot(S)$ must be a limit ordinal; hence by requirement \ref{MorCondBij},
$\tla>1$, and so by M.4, $\lambda$ is a limit ordinal.  
Also, from M.2, requirement \ref{MorCondBij},
and the fact that 
$\morl$ is a tree order, we get that for $x\morl y$ with $l(x)=\lambda$ and 
$l(y)=\ka$, $\pixy=f\restricted o(x)+1$.
Similarly, because $\morl$ is required to be a tree order and 
$f$ is surjective, one can only
have $x\morl y$ with $l(y)=\ka$ when 
either $l(x)=\lambda$, or there is a $z$ such that $x\morl z\morl y$ with
$l(z)=\lambda$.  
Thus, $\la$ is a mangal of the morass condition;  
this basic example of a mangal
will be important in later arguments.
As a result of $\la$ being a mangal,
it is equivalent to only explicitly require
M.3 to hold for those $y\in\calS$ with $l(y)<\ka$.
The same is trivially true of M.5, which for morass conditions
is vacuous in the $l(y)=\ka$ case.
Finally, note that because of the closure properties of $S$ given in
requirement \ref{MorCondLAS}, $f$ and $f^{-1}$ are SLOOP functions.

For any regular cardinal $\ka$,
let $P_\ka$ be the set of $(\ka,1)$-morass conditions, along with 
an extra point $\textbf{1}$ to
act as a maximum element; we will generally ignore 
$\textbf{1}$,
it being trivial to extend the definitions to encompass it as a special case.  
For each non-\textbf{1} element $p$ of $P_\ka$, let us denote with 
superscript $p$ the components and defined notions of $p$: 
$\calS^p$, $\morl^p$, $\pixy^p$ for $x\morl^py$, $\lambda^p$, 
$S^p$ as in requirement \ref{MorCondLAS}, 
$f^p$ as in requirement \ref{MorCondBij}, and $\tal^p$.
We shall also refer to
$\{x\in p\st l(x)=\lambda^p\}$
as the \emin{top level} of the condition $p$.

\begin{defn}\label{ManForcLEq}
\index{$\leq$ relation!mangrove forcing}
For $p$ and $q$ in $P$, we say that $q\leq p$ if
$$\calS^p\subseteq\calS^q,$$ 
$$\calS^p\cap((\lambda^p+1)\times\kappa)=
\calS^q\cap((\lambda^p+1)\times\kappa),$$
$$\morl^p = \morl^q\restricted\calS^p,$$
$$\forall x\morl^py\in\calS^p(\pixy^p=\pixy^q),$$
and $\la^p$ is a mangal of $q$.
\end{defn}
Notice in particular that a condition $p$ may be extended to a condition $q$
simply by extending $S^p$ and $f^p$ --- indeed, this is the only way to
extend $p$ without changing $\la$.  
This kind of extension is analogous to the possibility in 
Velleman's partial order of extending $p$ by increasing 
$\theta_{\lambda^p}$, but we avoid the necessity of 
(and indeed rule out) adding edges to the tree
below level $\lambda^p$.  We shall sometimes refer to the requirement 
that $\la^p$ be a mangal of $q$
(which is what is new in this forcing)
as the \emin{mangal requirement}; 
imposing it will ensure 
that $\la^p$ is a mangal of the generic morass if $p$ is in the generic,
as we shall see below.

The only slight difficulty in seeing that $\leq$ as defined above
is transitive is with the mangal requirement --- whether $\la^p$ is
a mangal of $r$ if $p\leq q\leq r$.  
This is easily resolved, however:

\begin{lemma}\label{ManLEqTrans}
The relation $\leq$ as given in Definition \ref{ManForcLEq} is transitive.
\end{lemma}
\begin{proof}
Suppose we have three conditions $p$, $q$ and $r$ in $P_\ka$ such that
$p\leq q\leq r$; we wish to show that $\la^p$ is a mangal of $r$.
So suppose there are $x$ and $y$ in $r$ such that $x\morl^ry$ and
$l(x)<\la^p<l(y)$.  If $l(y)\leq\la^q$, then $y\in q$ and so since 
$\la^p$ is a mangal of $q$, there is a $z\in q$ (and hence $z\in r$ also)
such that $l(z)=\la^p$ and $x\morl^q z\morl^q y$, whence
$x\morl^rz\morl^r$, as required.

If $l(y)>\la^q$, then because $\la^q$ is a mangal of $r$, there is some
$z'\in r$ such that $l(z')=\la^q$ and $x\morl^rz'\morl^ry$.  
But now we are reduced to the previous situation with $z'$ replacing $y$,
and so there is a $z\in r$ with $l(z)=\la^p$ and $x\morl^rz\morl^rz'$.
Transitivity of $\morl^r$ then gives $x\morl^rz\morl^ry$.
\end{proof}

Thus, letting $\P_\ka=\langle P_\ka, \leq\rangle$, we have that
$\P_\ka$\index{$\P_\ka$ (mangrove forcing)} is a partial order --- the
\emph{Mangrove Partial Order}\index{mangrove!partial order}
\index{partial order!Mangrove} or
\emph{Mangrove Forcing} 
at $\ka$.
In this section (Section \ref{MangroveForcing})
we will generally omit the subscript $\ka$, as $\ka$
will be unvarying.
Moreover, when we say that morass conditions are compatible or incompatible,
or comparable or incomparable, 
we will mean with respect to this mangrove partial order, rather than
Velleman's morass partial order.
We will also abuse notation, writing $p\in\P$ for $p\in P_\ka$.
The goal of the section, 
realised in Theorem \ref{limGMangrove},
is to show that forcing with $\P$ yields a mangrove in the generic extension.

\subsection{A first ``bamboo construction''}

It is not yet clear that there are any elements of
$\P$ other than $\boldsymbol{1}$.
The following lemma constructs one, and the 
basic technique of the construction, modified from one given in
\cite{Vell:MDF}, will be reused multiple times as we wish
to extend conditions in various ways later on.
We call such a construction a ``bamboo\index{bamboo} construction'' 
because of the
straight, vertical, unbranching branches that it produces.

\begin{lemma} \label{NonTrivMorCond}
There is a morass condition $q$ with $S^q=\omega$.
\end{lemma}
\begin{proof}
We start by defining $\tal^q$ for $\al$ up to $\omega^\omega$ 
(where the exponentiation is as ordinals), which will be our $\la^q$.
Let
$$
\tal=\begin{cases}
\{n\st\exists\zeta(\al=\omega^n\cdot\zeta)\}&
\text{if $\al<\omega^\omega$}\\
\omega&\text{if $\al=\omega^\omega$}\\
\end{cases}
$$

Now let
$$
\calS=
(\{\ka\}\times\omega)\cup\bigcup_{\al\leq\omega^\omega}(\{\al\}\times\tal),
$$
$$
x\morl^qy\longleftrightarrow l(x)<l(y)\land o(x)=o(y),
$$
and for $x\morl^qy$ and all $\ga\leq o(x)$,
$$
\pixy^q(\ga)=\ga.
$$
We claim that $q$ so defined is a morass condition.
The set $\calS$ clearly has the right form, $\morl^q$ is a (non-branching!) 
tree order, and for $x\morl^qy$, $\pixy^q$ is a function with the right domain
and range.  Requirements \ref{MorCondLAS} and \ref{MorCondBij} for a 
morass condition are immediate from the construction, as are
Monotonicity,
Commutativity,
M.1,
M.2,
M.5,
M.6 and M.7, the latter two being vacuous in the present context.
Finally,
axioms M.3 and M.4 below level $\ka$
reduce to simple properties of ordinal arithmetic:
respectively, 
that the limit of ordinals divisible by $\omega^n$ is divisible by
$\omega^n$; and that
if $\ga$ is divisible by $\omega^m$ for some $m>n$, then the set of ordinals in
$\ga$ divisible by $\omega^n$ is unbounded in $\ga$.
\end{proof}

The reader should note that such a bamboo construction does
not reflect the frequent branching that must occur in a mangrove to get
$\ka^+$-many leaves at level $\ka$.  
However, while one might not expect to find bamboo
in a mangrove, it does make a convenient building material.

\subsection{$\mu$-equivalence}\label{muequiv}

To acquaint the reader better with morass conditions and
mangrove forcing, we now present some
lemmas that may help with intuition and which will be useful later on.
First, we define a ``sameness'' notion which appears frequently,
and which will be crucial for the homogeneity arguments of Section 
\ref{homogeneity}.

\begin{defn}
Two morass conditions $p$ and $q$ are
\emph{$\mu$-equivalent},\index{equivalent!$\mu$- (morass conditions)}
\mbox{$p\mcequiv{\mu}q$}\index{$\mcequiv{\mu}$}, 
if they are ``the same up to level $\mu$'', that is,
\begin{eqnarray*}
\calS^p\cap((\mu+1)\times\ka)&=&\calS^q\cap((\mu+1)\times\ka),\\
\morl^p\cap((\mu+1)\times\ka)^2&=&\morl^q\cap((\mu+1)\times\ka)^2,
\end{eqnarray*}
for $x\morl^py$ with $l(y)\leq\mu$, 
$$
\pixy^p=\pixy^q,
$$
and $\mu$ is a mangal of both $p$ and $q$.
\end{defn}
Note that we include level $\mu$ in ``up to level $\mu$.''
Once again we shall refer to the last requirement in this definition as
the mangal requirement.

Clearly
$\mcequiv{\mu}$ is an equivalence relation on the set of those conditions with
$\la\geq\mu$ and $\mu$ as a mangal (although it is not even reflexive for other
conditions). 
While the presence of the mangal requirement means that 
$\mu\leq\nu$ does not yield
$p\mcequiv{\nu}q\implies p\mcequiv{\mu}q$, 
we  do have the following.

\begin{lemma}\label{mcequivIncrStr}
If $\mu\leq\nu$ and $p\mcequiv{\mu}q$ and $q\mcequiv{\nu}r$
for some elements $p$, $q$, $r$ of $\P$, then
$p\mcequiv{\mu}r$.
\end{lemma}
\begin{proof}
This is essentially the same as Lemma \ref{ManLEqTrans}, showing that
the $\leq$ relation of $\P$ is transitive.
Again, the only difficulty is in showing that $\mu$ is a mangal of $r$.
Suppose $x, y\in r$ are such that $x\morl^ry$ and
$l(x)<\mu<l(y)$. If $l(y)\leq\nu$, then $y\in q$ with 
$x\morl^qy$.  Thus, since $\mu$ is a mangal of $q$, there is a
$z\in q$ with $l(z)=\mu$ and $x\morl^qz\morl^qy$. So $z\in r$ and
$x\morl^rz\morl^ry$.

If $l(y)>\nu$, then since $\nu$ is a mangal of $r$, there is a
$z\in r$ with $l(z)=\nu$ and $x\morl^rz\morl^ry$.  
By transitivity of $\morl^r$, we are reduced to the case above.
\end{proof}

Also, it is immediate from the definitions that 
$q\leq p\implies q\mcequiv{\la^p}p$.
A stronger statement is in fact true.

\begin{lemma}\label{MorCondCompatMangal}
Let $p, q\in\P$ be (mangrove-) compatible morass conditions, 
with $\la^p\leq\la^q$.
Then $p\mcequiv{\la^p}q$.  
\end{lemma}
\begin{proof}
Let $r\in\P$ be such that $r\leq p$ and $r\leq q$.  
It is immediate from the fact that
$p\mcequiv{\la^p}r$ and
$q\mcequiv{\la^q}r$ that the necessary equalities for
$p\mcequiv{\la^p}q$ hold.  
To see that $\la^p$ is a mangal of $q$, 
observe that since $\la^p$ is a mangal of $r$,
we have that
for any $x$ and $y$ in $q$ with $l(x)<\la^p<l(y)$,
\begin{eqnarray*}
x\morl^qy&\longleftrightarrow& x\morl^ry\\
&\longleftrightarrow& \exists z\in r(l(z)=\la^p\land x\morl^rz\morl^ry)\\
&\longleftrightarrow& \exists z\in q(l(z)=\la^p\land x\morl^qz\morl^qy),
\end{eqnarray*}
so $\la^p$ is also a mangal of $q$.
\end{proof}
Note that, in the situation of the lemma, the only possible
obstruction to $q\leq p$ holding is that $f^p$ might not factor through $f^q$.
Indeed, $S^p$ need not be a subset of $S^q$, in which case $f^p$ certainly
cannot factor through $f^q$.

One reason that
the notion of $\mu$-equivalent is important, particularly in the
case where $\mu=\la^p$ for some condition $p$, is that most of the 
``morass-like'' requirements for the condition $p$ only pertain to
the part of the condition below $\la^p$.
The following lemma makes this observation concrete.

\begin{lemma}\label{MorCondChangeS}
Let $p$ be a morass condition, and let $S'$ be a subset of $\ka^+$
containing $0$, closed under ordinal successors and predecessors, and
such that $\ot(S')\leq\theta_{\la^p}^p$.  Then there is a unique
morass condition
$q$ such that $\la^q=\la^p$, $q\mcequiv{\la^p}p$, and $S^q=S'$.
\end{lemma}
\begin{proof}
The description of $q$ in the statement of the lemma completely 
determines $q$: to get from $p$ to $q$ we simply change $S^p$ to $S'$,
correspondingly replace $f^p$ with the order-preserving bijection
$f^q:\ot(S')\to S'$, and modify $\morl$ and $\pi$ as appropriate for 
the new $f$.  To be precise,
$$
\calS^q=\big(\calS^p\cap((\la^p+1)\times\ka)\big)\cup
(\{\ka\}\times S'),
$$
\begin{eqnarray*}
x\morl^qy&\Iff&\big(l(x)<l(y)\leq\la^p\ \land\ x\morl^py\big)\lor\\
&&\big(l(x)=\la^p\land l(y)=\ka\ \land\ f^q(o(x))=o(y)\big)\lor\\
&&\big(l(x)<\la^p\land l(y)=\ka\ \land\ 
\exists z\in\calS^q(l(z)=\la^p\land x\morl^qz\morl^qy)\big),
\end{eqnarray*}
and for $x, y\in\calS^q$ such that $x\morl^qy$,
$$
\pixy^q=\begin{cases}
\pixy^p&\text{if }l(x)<l(y)\leq\la^p\\
f^q\restricted(o(x)+1)&\text{if }l(x)=\la^p\text{ and }l(y)=\ka\\
\pi_{zy}^q\circ\pi_{xz}^q&\text{if }l(x)<\la^p, l(y)=\ka, l(z)=\la^p
\text{ and }x\morl^qz\morl^qy.
\end{cases}
$$
Every aspect of this definition is clearly
necessary for $q$ to be a morass condition
of the desired form, giving the claimed uniqueness.
If we can further show that this definition is sufficient to make
$q$ a morass condition we will be done.
But now the restricted forms of M.3--M.7 that are required all hold because 
they do in $p$ (recall from the
discussion after Definition \ref{MorassCondition} that M.3 and M.5 need
only be checked up to level $\la$).
The rest of the requirements for a morass condition are almost all 
immediate from the definitions; the only points worth mentioning are 
that $f^q$ is a SLOOP function,
and that M.2 in the case that $l(x)<\la^p<l(y)=\ka$ follows from
Commutativity and and M.2 in the other two cases of the definition of 
$\pixy^q$.
\end{proof}
In general this $q$ and the original $p$ will not be compatible ---
if there is a $\tau<\ka^+$ and $\al\neq\be<\theta_{\la^p}^p$ such that
$f^p(\al)=\tau=f^q(\be)$, there can be no common extension $r$, as it
would need to satisfy both 
$\langle\la^p,\al\rangle\morl^r\langle\ka,\tau\rangle$ and
$\langle\la^p,\be\rangle\morl^r\langle\ka,\tau\rangle$,
violating the fact that $\morl^r$ is a tree relation and Monotonicity 
holds in $r$.  On the other hand, in 
Lemma \ref{ManForcDiffSCompat} below we shall consider a particular
circumstance in which
such $\tau, \al$ and $\be$ do not exist, and $p$ and $q$ can shown to be
compatible.

\subsection{Cardinal preservation}\label{MangroveForcCardPres}

Our claim is that forcing with $\P$ will yield a 
mangrove with height $\ka$, with of course $\theta_\ka=\ka^+$.  
For this it is
important that the cardinals $\ka$ and $\ka^+$ are preserved.  In fact, 
if we assume that we have the GCH in our original model, all 
cardinals will be preserved, thanks to closure and chain condition properties
of $\P$ which we describe here.

Recall that a partial order $P$ is \emph{$\ka$-closed}\index{closed}
if every descending chain of conditions from $P$
of length less than $\ka$ has a lower bound in $P$.
Further recall that if $P$ is $\ka$-closed, then forcing with $\P$ adds no
new sequences of length less than $\ka$, and so in particular all cardinals
$\leq\ka$ are preserved.
In our particular circumstance,
knowing that $\P$ is $\ka$-closed will also be useful for ``gluing together''
successive extensions, making it possible to show that a variety of subsets
of $\P$ are dense using just a few basic extension lemmas repeatedly.

One may show that $\P$ is $\ka$-closed by the
natural argument, taking the union of the given chain of conditions, 
and adding a
top level as appropriate for the new $\ka$-th level.  As pointed out in 
\cite{Vell:MDF} though, 
the same argument actually yields a stronger property for
$\P$, one which will be useful when showing later that morasses give rise to
mangroves.  Perhaps more significantly, it will be exactly what we 
require in order to obtain master conditions when we force to obtain
mangroves while preserving large cardinals.

We make a preliminary definition in order to give Proposition
\ref{MangroveForcingKClosed} the 
necessary strength for the argument of Section \ref{ManfrMor}.

\begin{defn}\label{Pal}
For any $\al\leq\ka^+$, $\P_{\ka,\al}$ denotes the suborder of $\P$
consisting of those conditions $p\in\P$ such that $S^p\subseteq\al$.
\end{defn}

One might wonder whether the inclusion
$\P_{\ka,\al}\hookrightarrow\P$ is a complete embedding, as 
defined for example in \cite{Kun:ST}, page 218.
It turns out that this depends on $\al$ --- see Proposition
\ref{ManForcComplEmbs} in Section \ref{AsideComplEmbs}.

Recall that a set $Y$ in a partial order is \emin{directed}
if for every $p$ and $q$ in $Y$, there is an $r$ in $Y$ such that
$r\leq p$ and $r\leq q$.  Recall further that a partial order $P$ is
\emph{$\ka$-directed-closed}\index{directed-closed, $\ka$-} 
\index{closed!directed-}
if every directed subset
$Y$ of $P$ with cardinality less than $\ka$ has a lower bound in $P$.
Note in particular that the $\ka$ in ``$\ka$-directed-closed''
refers to the level of closure rather than the level of
directedness.  Of course, $\ka$-directed-closure implies 
$\ka$-closure.

\begin{prop} \label{MangroveForcingKClosed}
For every ordinal $\al\leq\ka^+$,
the poset $\P_{\ka,\al}$ is $\ka$-directed-closed.
In particular, $\P=\P_{\ka,\ka^+}$ is $\ka$-closed, 
and forcing with $\P$ adds no new sequences of length
less than $\ka$ of elements of $V$.
\end{prop}
\begin{proof}
Suppose  
that $Y$ is a directed set in $\P_{\ka,\al}$ with 
$|Y|=\gamma$ for some cardinal
$\gamma<\kappa$. 
Let $\bar q$ be the direct limit of these structures --- that is, let 
$\calS^{\bar q}=\bigcup_{p\in Y}\calS^{p}$, 
$x\morl^{\bar q}y$ if and only if
for some $p\in Y$ containing both $x$ and $y$, 
$x\morl^{p}y$, and in this case $\pixy^{\bar q}=\pixy^{p}$.
Note that $\bar q$ is well defined:  if $p,r\in Y$ both contain $x$ and $y$,
then a common extension $s\in Y$ will also contain $x$ and $y$,
and by the definition of $\leq$ in $\P$, $p$ and $r$ must both agree with
$s$ regarding whether $x\morl y$, and on the function $\pixy$ when 
$x\morl y$ does hold.

Suppose that $\bar q$ has a non-empty top level $\la^{\bar q}$. 
Then there exists at least one condition $p$ in $Y$ such that 
$\la^p=\la^{\bar q}$.
Recall that the only possibility for extending a condition $p$ without
altering $\la$ is to end extend $S^p$ and $f^p$.  Since $Y$ is directed,
it follows that conditions $p$ and $r$ in $Y$ with 
$\la^p=\la^r=\la^{\bar q}$ must be comparable: an $s\in Y$ extending both
simply has an $S^s$ end-extended from each of $S^p$ and $S^r$, and so 
the larger of the latter two sets must simply end-extend the smaller.
Thus, 
the subset $Y_{\la^{\bar q}}$ of $Y$ consisting of all elements $p$ of $Y$
with $\la^p=\la^{\bar q}$ forms a decreasing chain.  
Moreover, every element of
$Y$ is extended by an element of $Y_{\la^{\bar q}}$, so $\bar q$ is the direct
limit of the chain $Y_{\la^{\bar q}}$.
In this case it is 
straightforward to verify that all of the requirements for a morass condition
hold for $\bar q$, simply because they do for each $p\in Y_{\la^{\bar q}}$.  
Similarly, it is immediate from the definitions that
$\bar q$ is an extension of each $p\in Y_{\la^{\bar q}}$, and hence of
each $p\in Y$.
Finally, since we have simply taken a union to obtain $S^{\bar q}$,
$S^{\bar q}\subseteq \al$ and $\bar q$ is 
a lower bound for $Y$ in $\P_{\ka,\al}$.

So suppose now that $\bar q$ does not have a top level; 
then it is not a morass condition and we must modify it to obtain one.
Letting 
$\lambda^q=\sup_{p\in Y}(\lambda^{p})$,
$S^q=\bigcup_{p\in Y}S^{p}=
\{x\in\calS^{\bar q}\st l(x)=\ka\}$, 
and $f^q$ be the order preserving bijection from $\ot(S^q)$ to $S^q$,
define $q$ by 
$$
\calS^q=\calS^{\bar q}\cup(\{\lambda^q\}\times\ot(S^q)),
$$
\begin{eqnarray*}
x\morl^qy&\longleftrightarrow&\Big(x\morl^{\bar q}y\Big)\lor\\
&&\Big(l(x)=\lambda^q\land l(y)=\ka\land f^q(o(x))=o(y)\Big)\lor\\
&&\Big(l(x)<\la^q\land l(y)=\lambda^q\land \\
&&\qquad\qquad
\exists z\big(l(z)=\ka\land x\morl^{\bar q}z\land f^q(o(y))=o(z)\big)\Big),
\end{eqnarray*}
and for $x\morl^qy$,
$$
\pixy^q=\begin{cases}
\pixy^{\bar q} & \text{if $x\morl^{\bar q}y$,}\\
f^q\restricted(o(x)+1)& \text{if $l(x)=\lambda^q$ and $l(y)=\ka$,}\\
(f^q)^{-1}\circ\pi_{xz}^{\bar q}&
\text{if $l(y)=\lambda^q$, $l(z)=\ka$ and $x\morl^{\bar q}z$.}
\end{cases}
$$

We claim that with this definition $q$ is a morass condition.  
It is clear that $\calS^q$ has the right form, $\morl^q$ is a tree order,
and each $\pixy^q$ has the right domain and range.
Requirement \ref{MorCondBij} for a morass condition is also immediate from the
definition of $q$.
Requirement \ref{MorCondLAS} follows from the fact
that it is true for each $p\in Y$ and the easy observation that it 
also holds at the new top level. 
Monotonicity, Commutativity and M.1 are similarly straightforward
to verify.

For M.2, the cases where both $l(x)$ and $l(y)$ are less than $\lambda^q$
or where $l(x)<\lambda^q$ and $l(y)=\ka$
follow from M.2 for the conditions $p\in Y$.
If $l(x)=\la^q$ and $l(y)=\ka$, the fact that M.2 is satisfied is immediate 
from the definition of $q$.
Finally, the case of $l(x)<\la^q$ and $l(y)=\la^q$ is straightforward from the
fact that M.2 is satisfied in some $p\in Y$ containing $x$ and $z$, where 
$x\morl^qy\morl^qz$:
if $\nu<o(x)$, then 
$$\langle l(x),\nu\rangle\morl^{p}\langle\ka,\pi_{xz}^{p}(\nu)\rangle$$
and so
$$\langle l(x),\nu\rangle\morl^q
\langle\la^q,(f^q)^{-1}\circ\pi_{xz}^{p}(\nu)\rangle=
\langle l(y),\pixy^q(\nu)\rangle
$$ 
with
$$\pi^q_{\langle l(x),\nu\rangle
\langle l(y),\pixy^q(\nu)\rangle}=
(f^q)^{-1}\circ\pi_{xz}^{p}\restricted(\nu+1)=\pixy^q\restricted(\nu+1)$$
as required.

In the verifications of both M.3 and M.4, 
the only case we need to consider is when $l(y)=\la^q$, 
since for any $\al<\la^q$, any condition $p\in Y$ with $\la^p\geq\al$ fixes
the structure of the tree below level $\al$. 
So suppose that
$l(y)=\la^q$, and that $y\morl^qz$ with $l(z)=\ka$.  
Then there is some $p_z\in Y$ such that
$z\in p_z$.  For every $r\in Y$ extending $p_z$,
\begin{eqnarray*}
\{l(x)\st x\morl^qy\}\cap(\la^{r}+1)&=&
\{l(x)\st x\morl^qz\}\cap(\la^{r}+1)\\
&=&
\{l(x)\st x\morl^{r}z\}\\
\end{eqnarray*}
which is closed in $\la^r+1$ by M.3 for $r$.
Moreover, 
$\la^{r}\in\{l(x)\st x\morl^{r}z\}$ for each $r\leq p_z$ in $Y$
since $l(z)=\ka$.
Because $\la^q=\sup_{p\in Y}(\la^{p})$ and $Y$ is directed, 
we have that the set
$\{l(x)\st x\morl^qy\}$ is closed and unbounded in $\la^q$, 
as desired for M.3 and M.4.  
Also, this implies that for $y$ with $l(y)=\la^q$,
it is never the case that $x\morl_i^qy$.  Hence, the restricted forms of M.6
and M.7 that we require are vacuous at level $\la^q$, and of course hold 
below that level because they do in each $p_\al$.

It therefore only remains to check M.5, and again,
it is sufficient to consider the case when $l(y)=\la^q$.
Suppose $\mu<o(y)$; we wish to show that
there is an $x\morl^qy$ and a $\nu<o(x)$ such that 
$\pixy^q(\nu)=\mu$.
Let $w=\langle\ka,f^q(\mu)\rangle$, 
let $z=\langle\ka,f^q(o(y))\rangle$ so that $y\morl^qz$, and let
$p\in Y$ be such that $w$ and $z$ are both in $p$.  
Let $v,x\in p$ with $l(v)=l(x)=\la^{p}$ be such that $x\morl^{p}z$
and $v\morl^{p}w$.
Then $o(v)<o(x)$, and 
$$\pi_{xz}^{p}(o(v))=f^{p}(o(v))=\pi_{vw}^{p}(o(v))
=o(w)=f^q(\mu).$$
Now by the definition of $\morl^q$ we have that $x\morl^qy$ and 
$\pixy^q(o(v))=
\mu$ as required.

Overall then, 
we have that $q$ is a morass condition.  
Moreover, we claim that
$q\leq p$ for every $p\in Y$.  Indeed, the only part of this claim
that is not 
immediate from the definition of $q$
is the mangal requirement, and even this is trivial to verify 
when $l(y)\neq\la^q$.  
So suppose $p\in Y$, and $x, y\in q$ with 
$x\morl^qy$,
$l(y)=\la^q$ and
$l(x)<\la^p$.
Let $w=\langle\ka,f^{q}(o(y))\rangle$; then 
$x\morl^qw$. 
Let $r\in Y$, $r\leq p$ be such that $w\in r$, and hence $x\morl^rw$.
By the definition of extension, $\la^p$ is a mangal of $r$, and so
there is a $z\in r$ with $l(z)=\la^p$ such that
$x\morl^rz\morl^rw$, and hence $x\morl^qz\morl^qy\morl^qw$.
Thus, $\la^p$ is a mangal of $q$ for each $p\in Y$, 
and so $q\leq p$ for each $p\in Y$.

Hence, we have constructed a lower bound $q$ for $Y$.
As before, $S^q$ is simply the
union of $\{S^p\st p\in Y\}$, so $S^q\subseteq \al$ and $q\in \P_{\ka,\al}$. 
Thus, $\P_{\ka,\al}$ is indeed $\ka$-directed-closed.
\end{proof}


Our proof that $\P$ satisfies the $\ka^+$-cc will require a more complicated
construction, building a lower bound for conditions which are the same below 
level $\ka$ but have different sets $S$ at level $\ka$.

\begin{lemma}\label{ManForcDiffSCompat}
Let $p$ and $q$ be morass conditions such that $\la^p=\la^q$
and $p\mcequiv{\la^p}q$. 
Suppose further that there is some $S^0\subset\ka^+$ such that both
$S^p$ and $S^q$ end-extend $S^0$, and 
$\min(S^q\smallsetminus S^0)\geq\sup(S^p)$.  
Then $p$ and $q$ are compatible in $\P$.
\end{lemma}
\begin{proof}
We shall build a common extension
$r$ of $p$ and $q$, using a construction 
similar to that in the proof of Lemma \ref{NonTrivMorCond}.

Let $S^r=S^p\cup S^q$ and let $f^r$ be the order preserving bijection from
$\ot(S^r)$ to $S^r$.
Let $\ga^0=\ot(S^0)$, $\ga^p=\ot(S^p)$ and let $\ga_q$ be such that
$\ot(S^q)=\ga^0+\ga_q$; with this notation, 
$\ot(S^r)=\ga^p+\ga_q$.
Note also that by requirements \ref{MorCondLAS} and \ref{MorCondBij} 
for morass conditions,
$\ga^0$, $\ga^p$ and $\ga_q$ are limit ordinals, and $\ga^p+\ga_q<\ka$.

Let $\la^r=\la^p+\omega^{\ga^p+\ga_q}$.  For $\al\leq\la^r$, define $\tal^r$ by
$$
\tal^r=\begin{cases}
\tal^p& \text{if $\al\leq\la^p$}\\
\{\eta\st\exists\zeta(\be=\omega^\eta\cdot\zeta)\}&
\text{if $\al=\la^p+\be$, $0<\be<\omega^{\ga^p+\ga_q}$}\\
\ga^p+\ga_q& \text{if $\al=\la^r$.}\\
\end{cases}
$$
Now, let
$$
\calS^r=(\{\ka\}\times S^r)\cup\bigcup_{\al\leq\la^r}(\{\al\}\times \tal),
$$
and define $\morl^r$ and the $\pixy^r$ 
(in each case defined only when $x\morl^ry$, 
even though we do not state that below)
according to the following cases.  Let $\id_\al$ denote the identity function
on $\al$.
\begin{itemize}
\item If $l(x)\geq l(y)$ then $x\notmorl^ry$

\item If $l(x)<l(y)\leq\la^p$, then 
$x\morl^ry\leftrightarrow
x\morl^py$, and $\pixy^r=\pixy^p$.

\item If 
$\la^p<l(x)<l(y)<\ka$, then 
$x\morl^ry\leftrightarrow
o(x)=o(y)$, and $\pixy^r=\id_{o(x)+1}$.

\item If
$\la^p=l(x)<l(y)<\ka$ then
\begin{itemize} 
\item if $o(y)<\ga^p$ then 
$x\morl^ry\leftrightarrow
o(x)=o(y)$ and $\pixy^r=\id_{o(x)+1}$, and
\item if $o(y)\geq\ga^p$ then
$x\morl^ry\leftrightarrow
\exists\zeta(o(x)=\ga^0+\zeta\land o(y)=\ga^p+\zeta)$ and
$$
\pixy^r(\al)=\begin{cases}
\al& \text{if $\al<\ga^0$}\\
\ga^p+\eta& \text{if $\al=\ga^0+\eta$.}\\
\end{cases}
$$
\end{itemize}

\item If
$l(x)<\la^p<l(y)<\ka$ then
$$
x\morl^ry\leftrightarrow
\exists z(l(z)=\la^p\land x\morl^rz\morl^ry),
$$ and
$\pixy^r=\pi_{zy}^r\circ\pi_{xz}^r$.

\item If
$\la^r=l(x)<l(y)=\ka$ then
$x\morl^ry\leftrightarrow
f^r(o(x))=o(y)$ and $\pixy^r=$\hbox{$f^r\restricted(o(x)+1)$.}

\item If
$l(x)<\la^r<l(y)=\ka$ then
$$
x\morl^ry\leftrightarrow
\exists z(l(z)=\la^r\land x\morl^rz\morl^ry),
$$ and
$\pixy^r=\pi_{zy}^r\circ\pi_{xz}$.

\end{itemize}

We must verify that $r$ so defined is indeed a morass condition.
That $\calS$ is of the right form, $\morl^r$ is a tree order and
$\pixy^r$ is defined where necessary with the right domain and range is clear.
Requirements \ref{MorCondLAS} and \ref{MorCondBij} for a morass condition 
as well as Monotonicity, Commutativity, M.5, and the restricted form of
M.6 required are
immediate from the construction; for M.6 note in particular that
$\langle\la^p,\ga^0\rangle$ is $\morl^r$-less than both
$\langle\la^p+\omega^{\ga^p},\ga^0\rangle$ and
$\langle\la^p+\omega^{\ga^p},\ga^p\rangle$.
Once we have observed that Commutativity holds, M.2 also follows, as does
M.1 with the earlier observation that $\ga^p$ and $\ga^0$ are limit ordinals.
Axioms M.3 and M.4 are satisfied as in Lemma \ref{NonTrivMorCond} along with
the fact that they hold in $p$ and hence up to level $\la^p$ in $r$.
Finally, M.7 reduces to a property of ordinal arithmetic:
if $\de_0$ is a limit ordinal, and $\al$ is divisible by 
$\omega^\de$ for all $\de<\de_0$, then $\al$ is divisible by $\omega^{\de_0}$.

Therefore, $r$ is a morass condition.  
Clearly $r\leq p$ and $r\leq q$; 
in particular note that the mangal requirement holds
by construction.  Thus, $p$ and $q$ are compatible.
\end{proof}

\begin{prop}\label{MorForcK+cc} (GCH) The poset $\P$ is $\ka^+$-cc.
\end{prop}
\begin{proof} 
Suppose for the sake of contradiction that 
$A$ is an antichain in $\P$ of cardinality $\ka^+$.
Since $\la^p<\ka$ for any $p\in\P$, we may assume without loss of generality
that $\la^p=\la$ for all $p\in A$ and some $\la<\ka$.  Similarly, there are
at most $\ka^{|\la|}=\ka$ possible functions
$\theta:(\la+1)\rightarrow\ka$ determining the set $\calS$ up to level $\la$,
so by further paring down $A$ we may assume that 
$\calS^p\cap((\lambda+1)\times\ka)=
\calS^q\cap((\lambda+1)\times\ka)$
for all $p,q\in A$.
Furthermore, this common
initial segment of $\calS$ will have a cardinality $\mu<\ka$,
so there are at most $\ka$ many binary relations on it, and we may assume that
$\morl^p=\morl^q$ for $p,q\in A$.  With $\morl$ fixed, there are at most
$(\mu^\mu)^{\mu\times\mu}\leq\ka$ 
different possibilities for the set of functions
$\pixy$ with $l(x)<l(y)\leq\la$.  Hence, we may assume without loss of
generality that the conditions in $A$ are identical up to level $\la$.
Note that this implies that they are completely determined by the set
$S$ at level $\ka$: $S$ determines $f$, which in turn determines the relation
$x\morl y$ and maps $\pixy$ for $l(y)=\ka$.

It follows from the $\Delta$-system Lemma 
(see, for example, \cite{Kun:ST}, Theorem II.1.6)
that there is a $B\subset A$ of cardinality $\ka^+$ such that the set 
\hbox{$\{S^p\st p\in B\}$} forms a $\Delta$-system.
Moreover, $B$ may be chosen in such a way that there is a
$\ga<\ka$ and an $S^0\subset\ka^+$ such that for all $p, q\in B$,
$ 
\ot(S^p)=\ot(S^q)=\ga,
$ 
$S^p$ and $S^q$ end-extend $S^0$ with
$ 
S^p\cap S^q=S^0,
$ 
and 
$ 
\sup(S^p)<\min(S^q\smallsetminus S^0)$ 
or $\sup(S^q)<\min(S^p\smallsetminus S^0).
$ 
With $B$ so defined, any two of its elements will
be compatible by Lemma \ref{ManForcDiffSCompat}.  
This contradicts the initial assumption that $A$ was an antichain, and so we 
may conclude that $\P$ enjoys the $\ka^+$-cc.
\end{proof}

Note that we have in fact proved the slightly stronger property, that
$\P$ is \index{Knaster}\emph{$\ka^+$-Knaster}:
for any $A\subset\P$ of cardinality $\ka^+$, there is a $B\subset A$ of
cardinality $\ka^+$ whose elements are pairwise compatible.

\subsection{A mangrove from a generic}

We now turn to the main task of showing that the direct limit of the 
conditions in a $\P$-generic $G$ is a mangrove in $V[G]$.  
Of course, for this we will
need to show that various subsets of $\P$ are dense.  The following extension 
lemmas, in the same vein as Lemma \ref{NonTrivMorCond} and the argument for
Lemma \ref{ManForcDiffSCompat}, will give us the requisite density.

\begin{lemma}\label{DsigmaDense}  
For any $p\neq\boldsymbol{1}\in\P$ and any limit ordinal $\sigma<\ka^+$, there
is a $q\leq p$ such that $S^q=S^p\cup((\sigma+\omega)\smallsetminus\sigma)$.
\end{lemma}

\begin{proof}
If $\sigma\in S^p$, then we may of course take $q=p$, so suppose 
$\sigma\notin S^p$.  Let $S^q$ be as in the statement of the lemma.
If $\sigma\geq\sup(S^p)$ and 
$\ot(S^p)+\omega\leq\theta_{\lambda^p}$, we may simply extend $p$ to $q$ by
end-extending $f^p$, making the bijection from $\ot(S^p)+\omega$ to
$S^q$ our $f^q$.  Extending as necessary both the tree relation and the 
maps $\pixy$ by composition, it is straightforward to check that this defines
a $q\leq p$ with the desired $S^q$.  However, the construction described below 
works in all cases, so it is not necessary to distinguish this
particular case.

If the set $S^p\smallsetminus\sigma$ is non-empty, let $\tau$ be its least 
element, and let $\bar\tau$ be the unique ordinal such that 
$\langle\lambda^p,\bar\tau\rangle\morl^p\langle\ka,\tau\rangle$. 
Otherwise, let $\bar\tau=\ot(S^p)$.  
Clearly $\tau$ (if defined) and $\bar\tau$ are limit ordinals.
Further, let $\xi$ be such that $\bar\tau+\xi=\ot(S^p)$;
then $\xi$ is either $0$ or a limit ordinal, and 
$\ot(S^q)=\bar\tau+\omega+\xi$.

Note that $\bar\tau+\omega+\xi<\ka$, since 
$\bar\tau+\xi=\ot(S^p)\leq\theta_{\lambda^p}<\ka$, and so letting
$\delta=\omega^{\bar\tau+\omega+\xi}$, we also have $\delta<\ka$.  
Let $\lambda^q=\lambda^p+\delta$, and let
$$
\tal^q=\begin{cases}
\tal^p&\text{if $\al\leq\la^p$}\\
\{\nu \st \exists\zeta(\be=\omega^\nu\cdot\zeta\}&
\text{if $\al=\la^p+\be$, $0<\be<\delta$}\\
\bar\tau+\omega+\xi&\text{if $\al=\la^q$.}\\
\end{cases}
$$
Now let
$$\calS^q=\big(\bigcup\{\{\alpha\}\times\tal^q \st \al\leq\lambda^q\}\big)
\cup (\{\ka\}\times S^q);$$
then the usual definitions for $\tal^q$ and $S^q$ from $\calS^q$ match the 
chosen values, and we also immediately have
$f^q$, defined on the whole width $\theta_{\lambda^q}^q$ of the new 
top level $\lambda^q$.

Let us now define $\morl^q$ and the functions $\pixy^q$. 
This will be very similar to the definition for the construction in the
argument for Lemma \ref{ManForcDiffSCompat}.
%
As before, let $\id_\al$ denote the identity function
on $\al$, and although it is not explicitly stated, in each case take
$\pixy^q$ to be defined only when $x\morl^qy$.
\begin{itemize}
\item If $l(x)\geq l(y)$ then $x\notmorl^qy$

\item If $l(x)<l(y)\leq\la^p$, then 
$x\morl^qy\leftrightarrow
x\morl^py$, and $\pixy^q=\pixy^p$.

\item If 
$\la^p<l(x)<l(y)<\ka$, then 
$x\morl^qy\leftrightarrow
o(x)=o(y)$, and $\pixy^q=\id_{o(x)+1}$.

\item If
$\la^p=l(x)<l(y)<\ka$ then
\begin{itemize} 
\item if $o(y)\leq\bar{\tau}$ then 
$x\morl^qy\leftrightarrow
o(x)=o(y)$ and $\pixy^q=\id_{o(x)+1}$, and
\item if $o(y)>\bar{\tau}$ then
$x\morl^qy\leftrightarrow
\exists\nu(o(x)=\bar{\tau}+\nu\land o(y)=\bar{\tau}+\omega+\nu)$ and
$$
\pixy^q(\al)=\begin{cases}
\al& \text{if $\al<\bar{\tau}$}\\
\bar{\tau}+\omega+\mu& \text{if $\al=\bar{\tau}+\mu$.}\\
\end{cases}
$$
\end{itemize}

\item If
$l(x)<\la^p<l(y)<\ka$ then
$$
x\morl^qy\leftrightarrow
\exists z(l(z)=\la^p\land x\morl^qz\morl^qy),
$$ and
$\pixy^q=\pi_{zy}^q\circ\pi_{xz}^q$.

\item If
$\la^q=l(x)<l(y)=\ka$ then
$x\morl^qy\leftrightarrow
f^q(o(x))=o(y)$ and $\pixy^q=f^q\st(o(x)+1)$.

\item If
$l(x)<\la^q<l(y)=\ka$ then
$$
x\morl^qy\leftrightarrow
\exists z(l(z)=\la^q\land x\morl^qz\morl^qy),
$$ and
$\pixy^q=\pi_{zy}^q\circ\pi_{xz}$.

\end{itemize}

That completes the definition of $q$.
It is now a simple matter to check that the axioms for a morass condition 
hold for $q$; indeed, it is an identical verification process to that for the
$r$ of Lemma \ref{ManForcDiffSCompat}.  


So $q$ is a condition with 
$S^q=S^p\cup((\sigma+\omega)\smallsetminus\sigma)$, and by construction
$q$ extends $p$. 
\end{proof}

\begin{lemma} \label{MorCondVertExtn}
For any $p\neq\boldsymbol{1}\in\P$ there is a $q\leq p$ such that $S^q=S^p$,
$\theta_{\la^q}^q=\ot(S^q)$, and
$\la^q=\la^p+\omega^{\ot(S^q)}$.
\end{lemma}

\begin{proof}
The proof is another variant of the same argument, this one
being particularly straightforward. 
Let $S^q$ and $\la^q$ be as in the statement, set
$$
\tal^q=\begin{cases}
\tal^p&\text{if $\al\leq\la^p$}\\
\{\nu \st \exists\zeta(\be=\omega^\nu\cdot\zeta\}&
\text{if $\al=\la^p+\be$, $0<\be<\omega^{\ot(S^q)}$}\\
\ot(S^q)&\text{if $\al=\la^q$,}\\
\end{cases}
$$
and let
$$\calS^q=\big(\bigcup\{\{\alpha\}\times\tal^q \st \al\leq\lambda^q\}\big)
\cup (\{\ka\}\times S^q).$$
As always, take $f^q$ to be the order preserving bijection from
$\ot(S^q)$ to $S^q$.
Define $x\morl^qy$, and for $x$ and $y$ such that $x\morl^qy$, $\pixy^q$, 
by the following cases.  Again, $\id_\al$ denotes the identity function on 
$\al$.

\begin{itemize}
\item If $l(x)\geq l(y)$ then $x\notmorl^qy$

\item If $l(x)<l(y)\leq\la^p$, then 
$x\morl^qy\leftrightarrow
x\morl^py$, and $\pixy^q=\pixy^p$.

\item If 
$\la^p\leq l(x)<l(y)<\ka$, then 
$x\morl^qy\leftrightarrow
o(x)=o(y)$, and $\pixy^q=\id_{o(x)+1}$.

\item If
$l(x)<\la^p<l(y)<\ka$ then
$$
x\morl^qy\leftrightarrow
\exists z(l(z)=\la^p\land x\morl^qz\morl^qy),
$$ and
$\pixy^q=\pi_{zy}^q\circ\pi_{xz}^q$.

\item If
$\la^q=l(x)<l(y)=\ka$ then
$x\morl^qy\leftrightarrow
f^q(o(x))=o(y)$ and $\pixy^q=f^q\st(o(x)+1)$.

\item If
$l(x)<\la^q<l(y)=\ka$ then
$$
x\morl^qy\leftrightarrow
\exists z(l(z)=\la^q\land x\morl^qz\morl^qy),
$$ and
$\pixy^q=\pi_{zy}^q\circ\pi_{xz}$.

\end{itemize}

It is routine to check that with this definition $q$ is a morass condition
and $q\leq p$.
\end{proof}

%

Now let $G$ be $\P$-generic over $V$.  Let 
$M^G=\langle\calS^G,\morl^G,\langle\pixy^G\rangle_{x\morl^Gy}\rangle$
be the direct limit of the system $G$ in $V[G]$.  
That is, $\calS^G$ is the union over
$p\in G$ of $\calS^p$, 
for $x$ and $y$ in $\calS^G$ we set $x\morl^Gy$ if and only if there is some
condition $p\in G$ containing $x$ and $y$ such that $x\morl^py$, 
and in this case
$\pixy^G=\pixy^p$. 

\begin{thm} \label{limGMangrove}
With $M^G$ defined as above, $M^G$ is a mangrove in $V[G]$.
\end{thm}
\begin{proof}
Clearly $\calS^G$ is a subset of $(\ka\times\ka)\cup(\{\ka\}\times\ka^+)$
and the maps $\pixy^G$ are defined when they ought to be with the right 
domain and range.
The relation $\morl^G$ is tree relation since 
$\morl^p$ is a tree relation for each condition $p$, 
and below
level $\ka$ extensions of a condition only end-extend the relation.
From Lemma \ref{DsigmaDense}, 
we have that 
$D_\tau=\{p\in\P\st\tau\in S^p\}$
is dense in $\P$ for each
$\tau<\ka^+$,
and so it follows that $G\cap D_\tau\neq\emptyset$ and hence $S^G=\ka^+$. 
By Lemma \ref{MorCondVertExtn} along with the $\ka$-closure property of $\P$,
it also follows that for any $\al<\ka$, it is dense 
for the condition $p$ to have $\lambda^p>\al$.
In particular,
$\tal^G$ is an ordinal greater than $0$ for every $\al\leq\ka$, 
and the full strength of the Left-alignment axiom is satisfied.  

Monotonicity, Commutativity, M.1 and M.2 hold in $M^G$ 
because they do in each condition.
Similarly,
to show that M.3 and M.4 hold in $M^G$, it suffices to show that they do for
$y\in M^G$ with level $\ka$.  
So suppose $y\in M^G$ and $l(y)=\ka$.  Given $\al<\ka$, we
may take a condition $p$ such that $y\in p$ and $\lambda^p>\al$.  Then there
is some $z\in p$ with $l(z)=\lambda^p$ such that $z\morl^p y$.  
Hence, the set
$\{l(z)\st z\morl^G y\}$ is unbounded in $\ka$, as required for M.4.
Also, since the corresponding set in each condition is closed in $\ka$, and
the top levels $\lambda^p$ of the conditions are unbounded in $\ka$, 
we also have M.3: that the set is closed in $\ka$.

Axiom M.5 is also immediate everywhere in $M^G$ except level $\ka$,
from the fact that it holds in each condition.  
So suppose $y\in M^G$ has $l(y)=\ka$, and let
$\al<o(y)$.  Let $w=\langle\ka,\al\rangle$, and let $p\in G$ be such that
$w, y\in p$.  Let $v,x\in p$ be such that $l(v)=l(x)=\la^p$, $v\morl^pw$ and
$x\morl^py$.  Then by the remarks following definition \ref{MorassCondition},
$\pixy^p(o(v))=\al$.  Thus, M.5 holds in $M^G$ at level $\ka$ as well.

Axioms M.6 and M.7 require more work to verify, but fortunately the situation
is the same as in Velleman's paper \cite{Vell:MDF}.  We follow the arguments
there, but recast them as transfinite
inductions on natural well orders of the
elements of $M^G$, rather than proofs by 
contradiction for the least counterexample.

To verify M.6, 
suppose that $x\morl^Gy\in M^G$ with $o(x)$ a limit ordinal, and that
M.6 holds for all such $x'\morl^Gy'$ with $y'$ lexicographically less than $y$;
that is, $l(y')<l(y)$ or $l(y')=l(y)\land o(y')<o(y)$.  Let
$\zeta=\sup(\pixy^G``o(x))$ and let $z=\langle l(y),\zeta\rangle$.  Clearly if
$z=y$ we are done, so suppose $\zeta<o(y)$.

If $y$ is a limit point of $\morl^G$, then M.5 holds of $y$ and
there is some $w\in M^G$ and 
$\mu<o(w)$ such that $x\morl^Gw\morl^Gy$ and $\pi_{wy}(\mu)=\zeta$.
Denote $\langle l(w),\mu\rangle$ by $u$; then by axiom M.2, 
$u\morl^Gz$, and $\pi_{uz}^G=$\hbox{$\pi_{wy}^G\restr(\mu+1)$}.
By the composition axiom and M.1, we see that 
$\pi_{xw}^G``o(x)\subset\mu$, and moreover since 
$\pi_{xy}^G``o(x)$ is unbounded in $\zeta=\pi_{wy}^G(\mu)$,
$\pi_{xw}^G``o(x)$ must be unbounded in $\mu$.  
But now by the induction hypothesis, $x\morl^Gu$ with
$\pi_{xu}^G\restr o(x)=\pi_{xw}^G\restr o(x)$.  Hence composing, we have that
$x\morl^Gz$ and 
$$
\pi_{xz}^G\restr o(x)=\pi_{uz}^G\circ\pi_{xu}^G\restr o(x)=
\pi_{wy}^G\restr (\mu+1)\circ\pi_{xw}^G\restr o(x)=
\pi_{xy}^G\restr o(x)
$$
as required.

On the other hand, if $y$ is not a limit point of $\morl^G$, we may take a
$w\in M^G$ such that $w\morl_i^Gy$.  
Note that by M4, $l(y)\neq\ka$ in this case.
If $w=x$ then M.6 for $x$ and $y$
follows from the fact that it holds in any condition in $G$ containing both
$x$ and $y$.  
So suppose that $x\neq w$ and hence $x\morl^Gw\morl_i^Gy$.
Let $\mu=\sup(\pi_{xw}^G``o(x))$ and $u=\langle l(w),\mu\rangle$, 
so that by the induction hypothesis, $x\morl^Gu$ with 
$\pi_{xu}^G\restr o(x)=\pi_{xw}^G\restr o(x)$.  Note that it is possible that
$u=w$.  Let $\nu=\pi_{wy}(\mu)$ and $v=\langle l(y),\nu\rangle$; then
if $\mu<o(w)$, $u\morl^Gv$ with $\pi_{uv}^G=$
\hbox{$\pi_{wy}^G\restr (\mu+1)$} by M.2,
and otherwise $u=w$, $v=y$ and the same statement trivially holds true.
But now since $\pi_{xw}^G``o(x)$ is by definition unbounded in $\mu$,
\begin{eqnarray*}
\sup(\pi_{uv}^G``\mu)&=&\sup(\pi_{uv}^G\circ\pi_{xw}^G``o(x))\\
&=&\sup(\pi_{wy}^G\restr (\mu+1)\circ\pi_{xw}^G``o(x))\\
&=&\sup(\pi_{xy}^G``o(x))\\
&=&\zeta
\end{eqnarray*}
If $\mu<o(w)$ and $\nu<o(y)$ we have by the induction hypothesis that
$u\morl^Gz$ and $\pi_{uz}^G\restr \mu=\pi_{uv}^G\restr \mu$.  
Otherwise, $u=w$ and
$v=y$, so that $u\morl_i^Gv$ and we also obtain $u\morl^Gz$ and 
$\pi_{uz}^G\restr \mu=\pi_{uv}^G\restr \mu$, 
this time from the fact that it holds in any
condition containing both $u$ and $v$.
So then $x\morl^Gz$, and
$$
\pi_{xz}^G\restr o(x)=\pi_{uz}^G\circ\pi_{xu}^G\restr o(x)
=\pi_{uv}^G\circ\pi_{xw}^G\restr o(x)
=\pi_{wy}^G\circ\pi_{xw}^G\restr o(x)
=\pixy^G\restr o(x).
$$
Thus, we have by transfinite induction on the
lexicographic order for $y$ that M.6 holds of $M^G$.

To verify M.7 for $M^G$ we will induct simply on $l(y)$.
Let $x, y$ and $\al$ be as in the antecedent of M.7, and suppose first that
there is some $w$ such that $x\morl^Gw\morl^Gy$ and $l(w)>\al$.  
This is in particular the case when $y$ is a limit node of the tree,
thanks to Monotonicity and M.3.
For each $\nu<o(x)$, we have from M.2 that 
$$
\langle l(x),\nu\rangle\morl^G
\langle l(w),\pi_{xw}^G(\nu)\rangle\morl^G
\langle l(y),\pi_{wy}^G\circ\pi_{xw}^G(\nu)\rangle=
\langle l(y),\pixy(\nu)\rangle.
$$
Thus, since $\pixy^G``o(x)$ is unbounded in $o(y)$, it must be the case that
$\pi_{xw}^G``o(x)$ is unbounded in $o(w)$.  But now $x, w$ and $\al$ fit the
hypotheses for M.7, and $l(w)<l(y)$, so by the induction hypothesis, there
is a $\ga$ such that $\langle\al,\ga\rangle\morl^Gw\morl^Gy$.

So now suppose we are in the other situation, where there is some $w$ 
with $l(w)\leq\al$ such that $w\morl_i^Gy$.  If $l(w)=\al$ we 
are clearly done.  
Also, if $w$ were $x$, then M.7 would hold for $x$ and $y$ since it does for
any condition containing them, and that would contradict the assumption that
$\al$ is as in the hypothesis for that axiom.
So suppose that $x\morl^Gw\morl_i^Gy$ with $l(w)<\al$;  
we will derive a similar contradiction in this case.
As in the above case, we have for each $\nu<o(x)$ that
$$
\langle l(x),\nu\rangle\morl^G
\langle l(w),\pi_{xw}^G(\nu)\rangle\morl^G
\langle l(y),\pixy^G(\nu)\rangle,
$$
and $\pi_{xw}^G``o(x)$ is unbounded in $o(w)$.  
Let us denote by $\ga_\nu$ that ordinal such that
$\langle\al,\ga_\nu\rangle\morl^G\langle l(y),\pixy(\nu)\rangle$,
as we are assuming exists.  For simplicity, let us also define notation for
$\langle l(x),\nu\rangle,
\langle l(w),\pi_{xw}^G(\nu)\rangle,$ and
$\langle l(y),\pixy^G(\nu)\rangle$ for any given $\nu<o(x)$: call them
$x_\nu$, $w_\nu$ and $y_\nu$ respectively.
Then our string of tree relations can be rewritten and
extended thus:
$$
x_\nu\morl^G
w_\nu\morl^G
\langle \al, \ga_\nu\rangle\morl^G
y_\nu.
$$
Now, given $\mu<o(w)$, take some $\nu<o(x)$ such that 
$\pi_{xw}^G(\nu)>\mu$.  By M.2, we have that
$$
\langle l(w),\mu\rangle\morl^G
\langle\al,\pi_{w_\nu\langle\al,\ga_\nu\rangle}^G(\mu)\rangle
\morl^G
\langle l(y),\pi_{w_\nu y_\nu}^G(\mu)\rangle,
$$
and the corresponding functions $\pi$ are the restriction of 
$\pi_{w_\nu\langle\al,\ga_\nu\rangle}^G$ to $\mu+1$ and 
the restriction of $\pi_{\langle\al,\ga_\nu\rangle y_\nu}^G$ to 
$\pi_{w_\nu\langle\al,\ga_\nu\rangle}^G(\mu)+1$.
We now have that the hypotheses for M.7 hold of $w$, $y$ and $\al$, and
$w\morl_i^Gy$, so the conclusion also obtains, since it does in any
condition containing $w$ and $y$.  Of course, an ordinal $\ga$ given by
M.7 for $w$, $y$, and $\al$ would also act as a witness to the truth of 
M.7 for $x$, $y$, and $\al$.  
But in fact, we have a contradiction (as before)
to the assumption that $w\morl_i^Gy$, and so the supposed situation of 
having $w$ satisfying
$x\morl^Gw\morl_i^Gy$ and $l(w)<\al$ cannot occur.
In any case, we may conclude that M.7 holds in $M^G$.

That completes the verification that $M^G$ is a morass.  To confirm that it is
a mangrove, note that for any $p\in G$, $\la^p$ is a mangal of $M^G$: this is
immediate from the mangal requirement for extension.  
Thus, since the $\la^p$ for
$p\in G$ are unbounded in $\ka$, $M^G$ is a mangrove, as claimed.
\end{proof}

\section{A mangrove from a morass}\label{ManfrMor}

In this section we deviate from our main programme
in order to answer a natural question that arose earlier ---
whether the existence of a morass implies the existence of a mangrove.
A positive solution is provided by the theorems of 
\cite{Vell:MDF}, in which Velleman proves the equivalence of the existence
of a morass with a certain forcing axiom.  We shall show that
this forcing axiom is sufficient
to generate a mangrove from the mangrove forcing $\P$ defined in the previous 
section.  Indeed, our argument will closely follow Velleman's argument that
the forcing axiom gives a morass from his partial order.

We start with some preliminary definitions from \cite{Vell:MDF}.  
Let $P$ be a partial order,
and let $\calD=\{D_\al\st\al<\ka^+\}$ be an indexed family of open dense sets
in $P$.  For each $p\in P$, define the \emin{realm} of $p$, $\rlm(p)$, 
to be the set $\{\al<\ka^+\st p\in D_\al\}$.
For each $\al<\ka^+$ let $P_\al=\{p\in P\st\rlm(p)\subseteq\al\}$,
and let $P^*=\bigcup_{\al<\ka}P_\al$, the set of conditions with realms
bounded below $\ka$.

\begin{defn}[Definition 1.1.2 of \cite{Vell:MDF}]\label{AlmostKIndisc}
A family $\calD=\{D_\al\st\al<\ka^+\}$ of subsets of $P$ is said to
be \emin{almost $\ka$-indiscernible} if
\begin{description}
\item[I.1] $P^*\neq\emptyset$, and for all $\al<\ka$,
$D_\al\cap P^*$ is open dense in $P^*$, and
\item[I.2] for all $\al<\ka$, $P_\al$ is $\ka$-directed-closed,
\end{description}
and there is a function $\sigma$ assigning to each SLOOP
function $f:\al\to\ga$ with $\al<\ka$ and $\ga<\ka^+$
a function $\sigma_f:P_\al\to P_\ga$ such that the following properties 
also hold for all such $f$.
\begin{description}
\item[I.3] $\sigma_f$ is order preserving.
That is, 
$$
\forall p,q\in P_\al(p\leq q\implies\sigma_f(p)\leq\sigma_f(q)).
$$
\item[I.4] For all $p\in P_\al$, $\rlm(\sigma_f(p))=f``\rlm(p)$.
\item[I.5] If $\ga<\ka$, and
for some $\be<\al$, $f\restricted\be=\id_\be$ and $f(\be)\geq\al$,
then for all $p\in P_\al$, $p$ and $\sigma_f(p)$ are compatible in $P^*$.
\item[I.6] If $f_1:\al_1\to\al_2$ and $f_2:\al_2\to\ga$ are strictly
order preserving with $\al_1,\al_2<\ka$ and $\ga<\ka^+$, then
$$
\sigma_{f_2\circ f_1}=\sigma_{f_2}\circ\sigma_{f_1}.
$$
\end{description}
\end{defn}

One of the main results of \cite{Vell:MDF} is the following.

\begin{thm}[Theorems 1.1.3 and 2.1.6 of \cite{Vell:MDF}]\label{MorIffForcAx}
The following are equivalent:
\begin{enumerate}
\item There exists a $(\ka, 1)$-morass.
\item For any partial order $P$ and 
$\calD$ an almost $\ka$-indiscernible family
of subsets of $P$, there is a filter $G$ in $P$ such that
for every $D_\al\in\calD$, $G\cap D_\al\neq\emptyset$.
\end{enumerate}
\end{thm}

Velleman further strengthens (2) above by making $G$ hit $\ka$-many additional
sets $E_\zeta$; the reader trying to reconcile the statements here and in
\cite{Vell:MDF} may take each $E_\zeta$ to be $P$.  In this way
statement (2) above can be seen to lie in strength between the formulations in
Theorems 1.1.3 and 2.1.6 of \cite{Vell:MDF}, which are shown there to 
actually be equivalent.

To show that the existence of a morass implies the existence of a mangrove,
we will show that the subsets of the mangrove forcing
$$
D_\al=\{p\in\P_\ka\st\al\in S^p\}
$$
for $\al<\ka^+$ form an almost $\ka$-indiscernible family.
The filter we then obtain from Theorem \ref{MorIffForcAx} will 
have a mangrove as its direct limit, using essentially the
same proof as for Theorem \ref{limGMangrove}, where the filter 
in question was $\P$-generic over $V$.

Observe that with with these sets $D_\al$, $\rlm(p)=S^p$ for all
$p\in\P_\ka$, and $P_\al$ as defined before 
Definition \ref{AlmostKIndisc}
is simply the set of all morass
conditions $p$ with $S^p\subseteq\al$. 
That is, $P_\al=\P_{\ka,\al}$ as in 
Definition \ref{Pal}; we shall stick with the $\P_{\ka,\al}$ notation here.
Also, recall that for any morass condition $p$, 
$|S^p|\leq\theta_{\la^p}<\ka$, so if $p\in\P_{\ka,\ka}$, then there is
some $\al<\ka$ such that $p\in\P_{\ka,\al}$.  
Thus, $P^*=\P_{\ka,\ka}$; again, we will stick with the notation from
Definition \ref{Pal}.

\begin{prop}\label{ManForcDalAlmInd}
The sequence $\calD=\langle D_\al\st\al<\ka^+\rangle$ as defined above
is an almost $\ka$-indiscernible family of subsets of the mangrove 
forcing $\P_\ka$.
\end{prop}
\begin{proof}
Lemma \ref{NonTrivMorCond} shows that $\P_{\ka,\ka}$ is non-empty,
and Lemma \ref{DsigmaDense} shows that for each $\al<\ka$, 
$D_\al\cap\P_{\ka,\ka}$ is dense, with openness being clear from the
definition of $\leq$ in $\P$.  Therefore, I.1 holds.  Also,
I.2 is just Proposition \ref{MangroveForcingKClosed} 
for $\al<\ka$.

We must now define $\sigma_f:\P_{\ka,\al}\to\P_{\ka,\ga}$ for $f$ a SLOOP 
function with domain $\al<\ka$ and codomain $\ga<\ka^+$.  
If $p\in\P_{\ka,\al}$, let $\sigma_f(p)$ be that $q\in\P_{\ka,\ga}$ 
such that $\la^q=\la^p$, $q\mcequiv{\la^p}p$, and
$S^q=f``S^p$, as given by Lemma \ref{MorCondChangeS}.
Of course, this definition makes I.4 and I.6 immediate.
Further, the antecedent of I.5 makes Lemma \ref{ManForcDiffSCompat}
applicable, and the common extension constructed in that lemma
clearly lies in $\P_{\ka,\ka}$, so I.5 also holds.

To check I.3, let $g$ be a SLOOP function, and suppose $p\leq q$.  
Clearly $\la^q$ remains a mangal of $\sigma_g(p)$, and we have
$\sigma_g(q)\mcequiv{\la^q}\sigma_g(p)$.
If we can now show that 
\begin{multline*}
\morl^{\sigma_g(q)}\restricted
\big((\{\la^q\}\times\theta_{\la^q}^{\sigma_g(q)})\times
(\{\ka\}\times S^{\sigma_g(q)})\big)\\
=\ 
\morl^{\sigma_g(p)}\restricted
\big((\{\la^q\}\times\theta_{\la^q}^{\sigma_g(q)})\times
(\{\ka\}\times S^{\sigma_g(q)})\big),
\end{multline*}
and that for $x\morl^{\sigma_g(q)}y\in\calS^{\sigma_g(q)}$ with 
$l(x)=\la^q$ and $l(y)=\ka$,
$\pixy^{\sigma_g(q)}=\pixy^{\sigma_g(p)}$,
then we will be done, by Commutativity and the fact that 
$\la^q=\la^{\sigma_g(q)}$ is a mangal of $\sigma_g(q)$.
Further, since $f^{\sigma_g(q)}$ is surjective onto $S^{\sigma_g(q)}$,
it is sufficient to check that 
$\morl^{\sigma_g(q)}\subseteq\,\morl^{\sigma_g(p)}$ 
and 
$\pixy^{\sigma_g(q)}=\,\pixy^{\sigma_g(p)}$ on the relevant domain.

So suppose that $x\morl^{\sigma_g(q)}y$ for some $x, y\in\sigma_g(q)$ with
$l(x)=\la^q$ and $l(y)=\ka$.  
Note that by the definition of $\sigma_g(q)$,
$f^{\sigma_g(q)}=g\circ f^q$, and likewise with $q$ replaced by $p$.
Now
$$
x\morl^q\langle\ka,f^q(o(x))\rangle,
$$
so
$$
x\morl^p\langle\ka,f^q(o(x))\rangle.
$$
Let $z$ be that $z\in p$ such that
$l(z)=\la^p$ and 
$$
x\morl^pz\morl^p\langle\ka,f^q(o(x))\rangle.
$$
Then
$$
x\morl^{\sigma_g(p)}z\morl^{\sigma_g(p)}\langle\ka,g\circ f^q(o(x))\rangle
=\langle\ka,f^{\sigma_g(q)}(o(x))\rangle=y
$$
as required.

Similarly,
\begin{eqnarray*}
\pixy^{\sigma_g(q)}&=&f^{\sigma_g(q)}\restricted(o(x)+1)\\
&=&g\circ f^q\restricted(o(x)+1)\\
&=&g\circ\pi^q_{xy}\\
&=&g\circ\pi^p_{xy}\\
&=&g\circ f^p\circ\pi_{xz}^p\qquad\text{for $z$ as above}\\
&=&f^{\sigma_g(p)}\circ\pi_{xz}^p\\
&=&f^{\sigma_g(p)}\circ\pi_{xz}^{\sigma_g(p)}\\
&=&\pixy^{\sigma_g(p)}
\end{eqnarray*}
and we are done.
%
\end{proof}

From Proposition \ref{ManForcDalAlmInd} and Theorem \ref{MorIffForcAx},
we obtain that if a morass exists at $\ka$, then there is a filter $G$
in $\P_\ka$ intersecting every $D_\al$.  

In the proof of Theorem \ref{limGMangrove}, the only dense sets $G$ was 
required to meet were the sets $D_\al$, along with the open dense set of
conditions $p$ with $\la^p>\ga$, for each $\ga<\ka$.  It turns out that
using the following argument from \cite{Vell:MDF} 
we can avoid directly requiring
$G$ to meet these latter sets.

\begin{lemma}\label{ManForcGenEnoughCof}
If $G$ is a filter over $\P$ such that $G\cap D_\al\neq\emptyset$ for
each $\al<\ka^+$, then $\{\la^p\st p\in G\}$ is cofinal in $\ka$.
\end{lemma}
\begin{proof}
It will be convenient to work with the direct limit $M^G$ of the
filter $G$, with components $\calS^G$, $\morl^G$ and $\pi^G$.
Observe first that the verification 
of M.5 at level $\ka$ in $M^G$
in the proof of Theorem \ref{limGMangrove}
only required that $G$ intersect every $D_\al$, and thus goes through in
our situation.  
We also have that in $M^G$, $\morl$ is a tree order and
Monotonicity holds, since these properties are true in every condition,
and below level $\ka$ the relation \mbox{$\closermorl^p$} 
is only end-extended in
moving to an extension of $p$.
In particular, we have that for every $y\in \calS^G$,
$\{x\st x\morl y\}$ contains at most one element from each level.
Consider the node $y=\langle\ka,\ka\rangle$ in $M^G$.  
For every $x\morl^Gy$, $|\pixy^G``o(x)|=|o(x)|<\ka$.  Thus, since
$\ka$ is regular, there must be $\ka$ distinct $x\morl^Gy$, and 
$\{l(x)\st x\morl^Gy\}$ is cofinal in $\ka$.  But now for each $x\in \calS^G$,
there is a $p\in G$ such that $x\in p$, and hence $l(x)\leq\la^p$.  
It follows that $\{\la^p\st p\in G\}$ is cofinal in $\ka$.
\end{proof}

\begin{thm}
There is a morass at $\ka$ if and only if there is a mangrove at $\ka$.
\end{thm}
\begin{proof}
Suppose there is a morass at $\ka$.
As mentioned above, by proposition 
\ref{ManForcDalAlmInd} and Theorem \ref{MorIffForcAx},
there is a filter $G$ in $\P_\ka$ intersecting every $D_\al$.  
By Lemma \ref{ManForcGenEnoughCof}, $\{\la^p\st p\in G\}$ is cofinal in
$\ka$.
Knowing this, the proof of Theorem \ref{limGMangrove} goes through for
$G$, giving that the direct limit $M^G$ is a mangrove.  
Of course, the
other direction is trivial.
\end{proof}

\subsection{Aside: Complete embeddings}\label{AsideComplEmbs}

Although it is not necessary for the preceding arguments,
the curious reader might wonder whether the inclusion maps
$\P_{\ka,\al}\hookrightarrow\P_\ka$ are complete embeddings.
This seemingly simple question actually takes quite some work to 
answer, so we here devote a subsection to its resolution.

Recall the following definition from \cite{Kun:ST}.
\begin{defn}
Let $P$ and $Q$ be partial orders.  A function $i:P\to Q$ is a 
\emin{complete embedding} if
\begin{enumerate}
\item $\forall p,p'\in P\,(p'\leq p\implies i(p')\leq i(p))$
\item $\forall p,p'\in P\,(p\perp p'\iff i(p)\perp i(p')$
\item \label{reduction}
$\forall q\in Q\,\exists p\in P\,\forall p'\in P\,
(p'\leq p\implies i(p')\compat q)$
\end{enumerate}
A condition $p$ as in \ref{reduction} 
is said to be a \emin{reduction} of $q$ to $P$.
\end{defn}

\begin{prop}\label{ManForcComplEmbs}
For any limit ordinal $\al\leq\ka^+$, the inclusion 
$\P_{\ka,\al}\hookrightarrow\P_\ka$ is a complete embedding if and only if
$\cf(\al)=\ka$.
\end{prop}

Note that because of the requirement on morass conditions that $S$ be 
closed under ordinal successors, $\P_{\ka,\al+n}=\P_{\ka,\al}$ for any
limit ordinal $\al$ and natural number $n$.  The proposition therefore
effectively covers all cases.

\begin{proof}
Suppose first that $\al<\ka$.  Let $q\in\P_\ka$ be a morass condition such that
$\ot(S^q)>\al$ 
(such a $q$ may be constructed using lemmas
\ref{NonTrivMorCond} and \ref{DsigmaDense} and Proposition 
\ref{MangroveForcingKClosed}).
Let $p\in\P_{\ka,\al}$ be arbitrary; we shall show that $p$ is not
a reduction of $q$ to $\P_{\ka,\al}$.

Since $p\in\P_{\ka,\al}$, $\ot(S^p)\leq\al$.  Note that in the constructions
for Lemma \ref{MorCondVertExtn} and Proposition \ref{MangroveForcingKClosed},
the extension $q$ obtained satisfies $\theta_{\la^q}=\ot(S^q)$.
Thus, applying these results, we have that we can extend $p$ to a 
$p'\in\P_{\ka,\al}$
with $\la^{p'}>\la^q$ and $\theta_{\la^{p'}}^{p'}\leq\al$.
We claim that $p'$ is incompatible with $q$: if there were an $r\in\P_\ka$
extending both $p'$ and $q$, then $\la^{p'}$ would be a mangal of $r$, 
with $\theta_{\la^{p'}}^r=\theta_{\la^{p'}}^{p'}\leq\al$.
But the edge $\langle\la^q,\al\rangle\morl^q\langle\ka,f^q(\al)\rangle$
of $q$, and hence of $r$, cannot factor through such a level $\la^{p'}$
--- as a consequence of M.1, $o(x)\leq o(y)$ for any $x\morl^ry$.
Thus, no such common extension $r$ can exist, and
$p$ is not a reduction of $q$ to $\P_{\ka,\al}$.
Therefore
the inclusion $\P_{\ka,\al}\hookrightarrow\P_\ka$ is not a 
complete embedding.

Suppose next that $\cf(\al)=\ka$, 
and let $q$ be any element of $\P\smallsetminus\P_{\ka,\al}$.
Let $S^q_\al=S^q\cap\al$; because $\ot(S^q_\al)\leq\ot(S^q)<\ka$,
$S^q_\al$ is bounded in $\al$, and we will be able to overcome the
problems seen in the $\al<\ka$ case.
Let $\tau=\sup(S^q_\al)$ and $\mu=\ot(S^q\smallsetminus S^q_\al)$.
Using Lemma \ref{MorCondChangeS}, let $r$ be the condition such
that $\la^r=\la^q$, $r\mcequiv{\la^q}q$, and
$S^{r}=S^q_\al\cup((\tau+\mu)\smallsetminus\tau)$.
We claim that $r$ is a reduction of $q$ to $P_{\ka,\al}$.

Suppose $p\leq r$ in $\P_{\ka,\al}$, and note in particular that this 
implies that $S^p\supseteq S^r$.  
Applying Lemma \ref{MorCondChangeS} again, we can obtain a
$p'\in\P_\ka$ with $\la^{p'}=\la^p$, $p'\mcequiv{\la^{p}}p$, and
$S^{p'}=
(S^p\cap\tau)\cup
(S^q\smallsetminus S^q_\al)$.
Applying the construction of Lemma \ref{ManForcDiffSCompat},
we obtain a condition $s$ extending both $p$ and $p'$, with
$S^s=S^p\cup(S^q\smallsetminus S^q_\al)$.  We claim that $s\leq q$.
Certainly $\calS^s$, and in particular $S^s$, is large enough.  
We also have that 
$$
q\mcequiv{\la^q}r\mcequiv{\la^r}p\mcequiv{\la^p}s
$$
so $q\mcequiv{\la^q}s$.  Since $\la^q$ is a mangal of $q$, it 
suffices now to check that 
$$
\morl^q\restricted
\big((\{\la^q\}\times\theta_{\la^q}^q)\times(\{\ka\}\times S^q)\big)=\ 
\morl^s\restricted
\big((\{\la^q\}\times\theta_{\la^q}^q)\times(\{\ka\}\times S^q)\big),
$$
and that for $x\morl^qy\in\calS^q$ with $l(x)=\la^q$ and $l(y)=\ka$,
$\pixy^q=\pixy^s$.

Now for $x$ and $y$ in $q$ with $l(x)=\la^q$ and $l(y)=\ka$, 
$x\morl^qy$ if and only if $f^q(o(x))=o(y)$, and 
$\pixy^q=f^q\restr(o(x)+1)$. 
It therefore suffices
to show that for such $x$ and $y$, $x\morl^sy$ whenever 
$f^q(o(x))=o(y)$, with 
$\pixy^s=f^q\restr(o(x)+1)$ 
--- the converse holds because $f^q$,
the order preserving bijection from $\ot(S^q)$ to $S^q$,  
is surjective onto $S^q$.
So suppose $x\morl^qy$, that is, $f^q(o(x))=o(y)$.
There are two cases to consider.
\begin{itemize}
\item If $o(x)<\ot(S^q_\al)$, then $x$ and $y$ are actually in $r$, and
we have $x\morl^{r}y$ with $\pixy^r=\pixy^q$.
Since $r\geq p\geq s$, we have $x\morl^sy$ with $\pixy^s=\pixy^q$.
\item If $\ot(S^q_\al)\leq o(x)<\ot(S^q)$, let $\be$ be such that
$o(x)=\ot(S^q_\al)+\be$.  
Then in $r$, 
$$
x\morl^{r}\langle\ka,\tau+\be\rangle,\text{ with }
\pi_{x\langle\ka,\tau+\be\rangle}^r(\ga)=\begin{cases}
f^q(\ga)&\text{if }\ga<\ot(S^q_\al)\\
\tau+\de&\text{if }\ga=\ot(S^q_\al)+\de,\de<\be,
\end{cases}
$$ 
and so the same also holds true with $r$ replaced everywhere by $p$.
We therefore have that
$$
x\morl^p\langle\la^p,\ot(S^p\cap\tau)+\be\rangle\morl^p
\langle\ka,\tau+\be\rangle,
$$  
and 
$$
\pi_{x\langle\la^p,\ot(S^p\cap\tau)+\be\rangle}^p=
(f^p)^{-1}\circ\pi_{x\langle\ka,\tau+\be\rangle}^p.
$$
Recalling that $p'\mcequiv{\la^p}p$ and considering $f^{p'}$, 
we may conclude that
$$
x\morl^{p'}\langle\la^p,\ot(S^p\cap\tau)+\be\rangle\morl^{p'}
\langle\ka,f^q(\ot(S^q_\al)+\be)\rangle=y,
$$
with
$$
\pixy^{p'}=f^{p'}\circ
(f^p)^{-1}\circ\pi_{x\langle\ka,\tau+\be\rangle}^p.
$$
If $\ga<\ot(S^q_\al)$, then 
$$\pixy^{p'}(\ga)=f^{p'}\circ(f^p)^{-1}(f^q(\ga))=f^q(\ga).$$
If $\ga=\ot(S^q_\al)+\de$, then
\begin{eqnarray*}
\pixy^{p'}(\ga)&=&f^{p'}\circ(f^p)^{-1}(\tau+\de)\\
&=&
f^{p'}(\ot(S^p\cap\tau)+\de)\\
&=&f^q(\ot(S^q_\al)+\de)\\
&=&f^q(\ga).
\end{eqnarray*}
Hence, we have $x\morl^{p'}y$ with 
$\pixy^{p'}=f^q\restr(o(x)+1)=\pixy^q$,
and therefore the same is true in $s$.
\end{itemize}
In either case $s\leq q$ as claimed, and so $r$ is indeed a reduction
of $q$ to $\P_{\ka,\al}$.

Of course there is more to complete embeddings than the existence of 
reductions.  The inclusion $\P_{\ka,\al}\into\P$ certainly preserves
the relation $\leq$;  we must also show that it preserves incompatibility.
So suppose $p$ and $q$ in $\P_{\ka,\al}$ are compatible in $\P$ with some 
common extension $r\in\P$.  Let $r'$ be the unique condition such that
$\la^{r'}=\la^r$, $r'\mcequiv{\la^r}r$, and $S^{r'}=S^r\cap\al$, 
as given by Lemma \ref{MorCondChangeS}.  Clearly $r'\in\P_{\ka,\al}$,
and it is straightforward to verify that because $p\leq r$ and $q\leq r$,
$p\leq r'$ and $q\leq r'$.  So the inclusion $\P_{\ka,\al}\into\P$ 
indeed preserves incompatibility, and hence is a complete embedding.

Finally, suppose that $\al>\ka$ but $\cf(\al)<\ka$.  
Let $q\in\P\smallsetminus\P_{\ka,\al}$ be such that 
$S^q\smallsetminus\al\neq\emptyset$ and
$S^q\cap\al$ is cofinal in $\al$ (using Lemmas \ref{NonTrivMorCond} and
\ref{DsigmaDense} and Proposition \ref{MangroveForcingKClosed} as usual).
We claim that $q$ has no reduction to $\P_{\ka,\al}$.  
Suppose for the sake of contradiction that $p$ is a reduction of $q$ to
$\P_{\ka,\al}$, and let $r$ be a common extension of $p$ and $q$ in $\P$.
As in the argument above for incompatibility preservation, let $r'$ be the 
condition obtained from $r$ but truncating $S^r$ to $S^r\cap\al$.  
We again obtain $r'\leq p$, so $r'$ is also a reduction of $q$ to $\P$.
Also, if $q'$ is likewise obtained from $q$ by truncating $S^q$ to
$S^q\cap\al$, then it is straightforward to check that $q'\geq r'$
(indeed, $q'\geq q\geq r$, so the verification for incompatibility
preservation also applies to this case).  In particular, 
$S^{r'}\supseteq S^q\cap\al$, and for
any $\langle\la^q,\de\rangle$ in $q$ with $\de<\ot(S^q\cap\al)$,
$\langle\la^q,\de\rangle\morl^{r'}\langle\ka,f^q(\de)\rangle$. 

Let $r''$ be obtained from $r'$ by applying Lemma \ref{MorCondVertExtn},
so that in particular $\theta_{\la^{r''}}^{r''}=\ot(S^{r''})$ and
$S^{r''}=S^{r'}$.
Since the set of $f^q(\de)$ for $\de<\ot(S^q\cap\al)$ is cofinal in 
$\al=\sup(S^{r'})$, it follows that the set 
$$
\big\{\zeta\st
\exists\de<\ot(S^q\cap\al)\big(
\langle\la^q,\de\rangle\morl^{r''}\langle\la^{r''},\zeta\rangle
\morl^{r''}\langle\ka,f^q(\de)\rangle\big)\big\}
$$
is unbounded in $\theta_{\la^{r''}}^{r''}$.

Now suppose that $s$ is a common extension of $q$ and $r''$ with 
$\la^s>\la^{r''}$ (using Lemma \ref{MorCondVertExtn} again if necessary).
Let $\be=\ot(S^q\cap\al)$; 
then there is some $\tau\geq\al$ such that 
$\langle\la^q,\be\rangle\morl^q\langle\ka,\tau\rangle$, and hence
$\langle\la^q,\be\rangle\morl^s\langle\ka,\tau\rangle$.
Since $\la^{r''}$ is a mangal of $s$, let $\ga$ be such that 
$\langle\la^q,\be\rangle\morl^s\langle\la^{r''},\ga\rangle\morl^s
\langle\ka,\tau\rangle$.
We must have 
$
\ga<\theta_{\la^{r''}}^s=\theta_{\la^{r''}}^{r''},
$
so there is some $\de<\ot(S^q\cap\al)$ and $\zeta>\ga$ such that 
$\langle\la^q,\de\rangle\morl^{r''}\langle\la^{r''},\zeta\rangle
\morl^{r''}\langle\ka,f^q(\de)\rangle$.
Now 
$$
\pi_{\langle\la^{r''},\ga\rangle\langle\ka,\tau\rangle}^s\circ
\pi_{\langle\la^{q},\be\rangle\langle\la^{r''},\ga\rangle}^s(\de)
=\pi_{\langle\la^{q},\be\rangle\langle\ka,\tau\rangle}^s(\de)
=\pi_{\langle\la^{q},\be\rangle\langle\ka,\tau\rangle}^q(\de)
=f^q(\de)
$$
So by axiom M.2, $\langle\la^{r''},
\pi_{\langle\la^{q},\be\rangle\langle\la^{r''},\ga\rangle}^s(\de)\rangle
\morl^s\langle\ka,f^q(\de)\rangle$.
But
$\pi_{\langle\la^{q},\be\rangle\langle\la^{r''},\ga\rangle}^s(\de)<\ga
<\zeta$, and we already have that 
$\langle\la^{r''},\zeta\rangle\morl^s\langle\ka,f^q(\de)\rangle$, 
giving a contradiction.  Thus, there cannot be a reduction $p$ of $q$ to
$\P_{\ka,\al}$, 
and hence the inclusion $\P_{\ka,\al}\into\P$ is not a complete
embedding.
\end{proof}

\section{Homogeneity}\label{homogeneity}

When we force to obtain a morass in the presence of various large cardinal 
axioms, we will sometimes be able to preserve the large cardinal property by
simply requiring that our generic be below some master condition.  Of
course, for each cardinal $\ka$ with the given large cardinal property,
there will be a different master condition.  Hence, to be able to preserve 
the large cardinal property for all of these cardinals at once, 
one would like to be able to impose each master 
condition ``after the fact'', 
substituting it into the ``bottom'' of any given generic to give a
new generic suited to the large cardinal at hand.
This is achieved with Theorem \ref{MangroveForcingHomog} of this section.

In some forcing posets, such as the standard ones to force 
$\diamondsuit$ and $\square$, compatibility implies comparability, and
further one of the conditions will be an initial segment of the other.  
When these initial segments may be freely interchanged giving new 
conditions, a great deal of homogeneity for the partial order results.

While we do not have these same properties in $\P$, 
our notion of $\mu$-equivalence is an appropriate substitute for
comparability which does follow from compatibility 
(Lemma \ref{MorCondCompatMangal}), with the part of the morass condition
below $\mu$ being a suitable ``initial segment'' for us to interchange.
To be precise, we will
use $\mu$-equivalence in defining
automorphisms on suitably chosen open dense suborders of $\P$,
with the automorphism being a sort of ``initial segment interchange''.  
We will
then show in Theorem \ref{MangroveForcingHomog} that such automorphisms
may be used to construct, within
an arbitrary generic extension $V[G]$,
$\P$-generics over $V$
containing any given condition.
The definition of these automorphisms is the crux of our argument for
Theorem \ref{MangroveForcingHomog}, and by extension, the large cardinal
preservation arguments in Sections \ref{MangroveLCs} and
\ref{Mangrove1Extend}.  Velleman's partial order of \cite{Vell:MDF},
with no mangal requirement for extension, lacks this crucial homogeneity.

First, we define the suborders of interest.

\begin{defn}\label{MangrovePalpha}
\index{$\P^\al$ (mangrove forcing above $\la^\al$)}
For any $\al<\ka$, let 
$$
\P^\al=\{p\in\P\st\la^p\geq\al\}.
$$
\end{defn}

\begin{prop}\label{MangrovePalphaDense}
For all $\al<\ka$, $\P^\al$ is open dense in $\P$.
\end{prop}
\begin{proof}
Density follows from 
Lemma \ref{MorCondVertExtn} and Proposition \ref{MangroveForcingKClosed},
and openness is immediate from the definition of $\leq$.
\end{proof}

Of course, this means that generics for $\P^\al$ are easily
translatable into generics for $\P$ and conversely, while preserving the
extension universe $V[G]$.

We will define our automorphisms by ``interchanging''
initial pieces of conditions.
\begin{defn}\label{Mangroveqpr}
Let $p, q$ and $r$ in $\P$ be three conditions
such that 
$\la^{q}\geq\la^{p}=\la^{r}$, 
$\theta_{\la^{p}}^{p}=
\theta_{\la^{r}}^{r}$,
and
$q\mcequiv{\la^{p}}p$. 
Define $q^{p}_{r}$\index{$q^{p}_{r}$} as follows:
let
\begin{eqnarray*}
\tal^{q^{p}_{r}}&=&\begin{cases} 
\tal^{r}&\text{if $\al\leq\la^{p}$}\\
\tal^q&\text{if $\al\geq\la^{p}$}\\
\end{cases}\\
S^{q^{p}_{r}}&=&S^q
\end{eqnarray*}
\begin{eqnarray*}
x\morl^{q^{p}_{r}}y&\longleftrightarrow
 &\big(l(x)<l(y)\leq\la^{p}\enskip\land\enskip x\morl^{r}y\big)\enskip\lor\\
&&\big(\la^{p}\leq l(x)<l(y)\enskip\land\enskip x\morl^{q}y\big)\enskip\lor\\
&&\big(l(x)<\la^{p}<l(y)\enskip\land\\
&&\qquad\qquad
\exists z\in\calS^q(l(z)=\la^{p}\land x\morl^{r}z\morl^qy)\big),
\end{eqnarray*}
and for 
$x\morl^{q^{p}_{r}}y$,
$$
\pixy^{q^{p}_{r}}=\begin{cases}
\pixy^{r}&\text{if $l(x)<l(y)\leq\la^{p}$}\\
\pixy^q&\text{if $\la^{p}\leq l(x)<l(y)$}\\
\pi_{zy}^q\circ\pi_{xz}^{r}&
\text{if 
$l(z)=\la^{p}\land x\morl^{r}z\morl^qy$.}\\
\end{cases}
$$
\end{defn}

Intuitively $q^{p}_{r}$ is just $q$ with $p$ replaced by $r$,
although $S^p$ and $S^r$ are irrelevant for the construction. 

\begin{prop}\label{Mangroveqprok}
For $p, q$ and $r$ as in Definition \ref{Mangroveqpr}, $q^p_r$ is
a morass condition and $q^p_r\mcequiv{\la^p}r$.
\end{prop}
\begin{proof}
It is straightforward to check
that all of the requirements for a 
morass condition hold for $q^{p}_{r}$, because they hold for $q$ and $r$.
In particular, note that M.2 holds for $l(x)<\la^{p}<l(y)$ because 
$\pixy^{q^{p}_{r}}$ is defined by composition for such $x$ and $y$.  
Also, because we only require M.6 and M.7 to hold for immediate 
$\morl$-successors in a morass condition, these two axioms are immediate from
the fact that they hold in $q$ and $r$.
It is immediate from the definition that
$\la^{p}$ is also a mangal of the new condition
$q^{p}_{r}$, and $\la^p=\la^r$ so it is also a mangal of $r$.
Since the requisite equalities trivially hold, it follows that 
$q^{p}_{r}\mcequiv{\la^{p}}r$.
\end{proof}

\begin{defn}\label{Mangrovephipr}
Let $p$ and $r$ in $\P$ be such that $\la^p=\la^r$ and 
$\theta^p_{\la^p}=\theta^r_{\la^r}$.
Then define $\varphi^{p,r}:\P^{\la^p}\to\P^{\la^p}$ by
$$
\varphi^{p,r}(q)=\begin{cases}
q^p_r&\text{if $q\mcequiv{\la^p}p$}\\
q^r_p&\text{if $q\mcequiv{\la^p}r$}\\
q&\text{otherwise.}
\end{cases}
$$
\end{defn}
Note first that this is well defined: if 
$p\mcequiv{\la^p}q\mcequiv{\la^p}r$,
then 
$p\mcequiv{\la^p}r$, and it is easy to see that $q^p_r=q^r_p=q$.

Because of the mangal requirement for the relation $\mcequiv{\la^p}$, 
$\varphi^{p,r}$ is self-inverse: ${\varphi^{p,r}}^2(q)=q$ for all 
$q\in\P^{\la^p}$.
Note that this is the point in the argument which necessitates the
mangal requirement, and concomitantly where
the definition of $\leq$ in Velleman's partial
order fails for our purposes.  
We could have defined $q^{p}_{r}$ without such a
requirement, and would still have obtained a morass condition, 
but the corresponding $\varphi$ would not have been surjective.
The range of such a $\varphi$ would not include any
condition that was ``the same as $p$ up to level $\la^p$'' but
which did not have $\la^p$ as a mangal, and likewise for $r$.
Since the set of such conditions is not codense in $\P^{\la^p}$, 
this is a real cause for concern, hence the introduction of mangals.

\begin{prop}
For conditions $p$ and $r$ satisfying the requirements of 
Definition \ref{Mangrovephipr}, 
the function $\varphi^{p,r}$ is an automorphism of the poset $\P^{\la^p}$.
\end{prop}
\begin{proof}
Since $\varphi^{p,r}$ is a self-inverse bijection, it only remains to show
that $\varphi^{p,r}$ respects $\leq$.
So suppose that $q\leq s\in\P^{\la^p}$.
Then 
$q\mcequiv{\la^s}s$ and $\la^{p}<\la^s$, so 
$q\mcequiv{\la^{p}}p\iff s\mcequiv{\la^{p}}p$ and
$q\mcequiv{\la^{p}}r\iff s\mcequiv{\la^{p}}r$ by 
Lemma \ref{mcequivIncrStr}
(note that this would not necessarily be the case if we were to use 
compatibility in place of the $\mcequiv{\la^{p}}$ relation).

For the sake of the argument at hand, we may therefore assume without
loss of generality that $q\mcequiv{\la^p}p\mcequiv{\la^p}s$, so that
$\varphi^{p,r}(q)=q^p_r$ and
$\varphi^{p,r}(s)=s^p_r$.
To see that $\la^q$ remains a mangal of $\varphi^{p,r}(s)$, note by 
Proposition \ref{Mangroveqprok} that $\la^p=\la^r$ is a mangal of $s^p_r$,
so any edge of $s^p_r$ factors as an edge in $r$ followed by an edge in $s$
from level $\la^p$, the latter factorising through level $\la^q$ since
$\la^q$ is a mangal of $s$.  The remainder of the verification that
$q^p_r\leq s^p_r$ is trivial from the definitions.
\end{proof}

The following lemma will also be useful.
\begin{lemma}\label{phiprSwapspr}
If $p$ and $r$ are conditions in $\P$ such that
$\la^p=\la^r$ and 
$\theta^p_{\la^p}=\theta^r_{\la^r}$, and further $S^p=S^r$, then
$\varphi^{p,r}(p)=r$ and
$\varphi^{p,r}(r)=p$.
\end{lemma}
\begin{proof}
Of course $\varphi^{p,r}$ is symmetric in $p$ and $r$, so it suffices to 
consider one case.
The set $S^r$ determines $f^r$ and hence also $\morl^r$ and $\pixy^r$ up from
level $\la^r$, and $\varphi^{p,r}(p)=p^p_r\mcequiv{\la^r}r$.  
\end{proof}

With these tools at our disposal, we are almost ready to prove the main
theorem of this section.
The idea of the proof of Theorem \ref{MangroveForcingHomog}
is to find a condition $r$ in $G$ large enough to 
``cover'' a given $p$, so that we can use $\varphi^{p,r}$ to
``patch in'' $p$ in place of $r$ in the generic.  
The following two lemmas allow us to find a suitable such $r$.


\begin{lemma} \label{DXDense}
For any subset $X$ of $\ka^+$ such that $|X|<\ka$, the set
$$
D_X=\{r\in\P\st X\subset S^r\}
$$ 
is open dense in $\P$.
\end{lemma}
\begin{proof} 
Density is obvious by combining Lemma \ref{DsigmaDense} with Proposition
\ref{MangroveForcingKClosed}.  Openness is immediate from the definition of
extension.
\end{proof}

\begin{lemma} \label{MorCondEmu}
For any ordinal $\mu<\ka$, the set
$$
E_\mu=
\{r\in\P\st\exists\al>\ot(S^r)(\lambda^r=\mu+\omega^{\ot(S^r)}\cdot\al)
\land \theta_{\la^r}^r=\ot(S^r)\}
$$ 
is dense in $\P$.
\end{lemma}
\begin{proof} 
First note that $\mu$ only affects the minimum value of $\la$ for elements of 
$E_\mu$,
and not its ``form'': 
it is a simple result of ordinal arithmetic, that
for any $\mu$, $\ga$ and $\al_0$, there is an $\al_1$ such that
$\mu+\omega^\ga\cdot\al_0=\omega^\ga\cdot\al_1$.
Now given any $p$ in $\P$, we may simply
use Lemma \ref{MorCondVertExtn} repeatedly, ``gluing together'' with
Proposition \ref{MangroveForcingKClosed} (noting that the argument there 
will not increase $S$), to give an 
extension in $E_\mu$.
\end{proof}

We are now ready to prove the main theorem of this section.
\begin{thm}[Homogeneity]\label{MangroveForcingHomog} 
Let $G$ be $\P$-generic over $V$ and let $p$ be a condition in 
$\P$.  Then there exists a $G'\subset V[G]$ such that $p\in G'$ and 
$G'$ is also $\P$-generic over $V$.  Moreover, $V[G']=V[G]$.
\end{thm}
\noindent Note that we write $G'\subset V[G]$ rather than $G'\in V[G]$ because 
we are not assuming that $\P$ is set sized from the perspective of 
the ground model $V$ ---
this will be relevant in Section \ref{Mangrove1Extend} below.
\begin{proof}
By Lemmas \ref{DXDense} and \ref{MorCondEmu}, 
there is some $r\in G\cap D_{S^p}\cap E_{\la^p}$.
We claim that $p$ may be extended to a condition $p'$ 
such that $S^{p'}=S^r$, $\la^{p'}=\la^r$, and
$\theta_{\la^{p'}}^{p'}=\theta_{\la^r}^r$. 
Note that in Lemma \ref{DsigmaDense}, the extension $q$ of $p$
constructed satisfies
$\la^q=\la^p+\omega^{\ot(S^q)}$.
This construction is to add a single new limit ordinal to $S$, whereas
we wish to add the whole of $S^r\smallsetminus S^p$.  
However, 
we may construct
a sequence of conditions $q_\al$, 
whereby if 
$\langle\sigma_\al\st\al<\be\rangle$ is the increasing enumeration of 
limit ordinals in $S^r\smallsetminus S^p$, $q_{\al+1}$ is obtained from 
$q_{\al}$ by the method of Lemma 
\ref{DsigmaDense} adding $\sigma_\al$ and its finitary successors, 
$q_0=p$, and $q_\la$ for limit ordinals $\la$ is
obtained by the method of Proposition \ref{MangroveForcingKClosed} to close off
with a new top level.  
Using Proposition \ref{MangroveForcingKClosed} one last time if
$\be$ is a limit ordinal and taking $q_{\be-1}$ otherwise, 
we get a lower bound $q$ for the sequence 
such that $S^q=S^r$ and
$$
\la^q=\la^p+\sum_{0<\al<\be}\omega^{\ot(S^{q_\al})}.
$$

Now, the summation on the right (excluding $\la^p$) is 
clearly less than or equal to
$\omega^{\ot(S^q)}\cdot\be$.
Further, $\be$ must be less than or equal to $\ot(S^q)$, so we have that 
$$
\la^q\leq\la^p+\omega^{\ot(S^q)}\cdot\ot(S^q).
$$
But now 
$$
\la^r=\la^p+\omega^{\ot(S^r)}\cdot\de
=\la^p+\omega^{\ot(S^q)}\cdot\de
\quad\text{for some $\de>\ot(S^q),$}
$$
so noting as in Lemma \ref{MorCondEmu} that adding 
$\omega^{\ot(S^q)}$ to an ordinal yields an ordinal divisible by 
$\omega^{\ot(S^q)}$,
we see that
we may extend $q$ further to a $p'$ with
$S^{p'}=S^r$ and $\la^{p'}=\la^r$, by simply 
applying Lemma \ref{MorCondVertExtn} and 
Proposition \ref{MangroveForcingKClosed} repeatedly. 
Also, since we must use Lemma \ref{MorCondVertExtn} at least once to get $p'$, 
we will further have that 
$\theta_{\la^{p'}}^{p'}=\ot(S^{p'})=\ot(S^r)=\theta_{\la^r}^r$, 
and so $p'$ has all of the desired attributes in common with $r$.

Now in $V[G]$, let $G^{\la^r}$ denote $G\cap\P^{\la^r}$,
let $G^{\la^r}_{p'}$ denote $\varphi^{p',r}``G^{\la^r}$, and let
$$
G'=\{q\in\P|\exists s\in G^{\la^r}_{p'}(s\leq q)\}.
$$
Because $\varphi^{p',r}$ is a poset automorphism and 
$\P^{\la^r}$ is open dense in $\P$, $G'$ is $\P$-generic over $V$, and
$V[G']=V[G]$.  Because $\varphi^{p',r}(r)=p'$ and
$r\in G^{\la^r}$, $p'\in G^{\la^r}_{p'}$, and so $p\in G'$,
as desired.
\end{proof}


\section{Preserving large cardinals}\label{MangroveLCs}

In \cite{SDF:LCL}, it was shown that using a reverse Easton iteration of
forcings to construct morasses, one may force morasses to exist at every 
cardinal while preserving a hyperstrong or $n$-superstrong cardinal, for 
$0<n\leq\omega$.  However, except in the case of 1-superstrong cardinals,
this result is 
only achieved for a single cardinal, as
the generic used must fall below a specific 
master condition that is dependent on the embedding witnessing 
$n$-superstrength or hyperstrength, and 
so a single generic cannot be used to preserve the large cardinal strength
of many different cardinals.  
On the other hand, with our Homogeneity Theorem
(Theorem \ref{MangroveForcingHomog}), 
we can circumvent this problem, expressing 
for each master condition $p$ the 
extension $V[G]$ as $V[G']$ for some $G'\ni p$.  
We give the details of the argument in this section.

Assume the GCH holds (this can be forced while preserving
large cardinals as in \cite{SDF:LCL} and Section \ref{GCH}).
\begin{defn}\label{GlobManPO}
The \emin{Global Mangrove Partial Order}\index{partial order!Global Mangrove} 
$R$ is the reverse Easton iteration
of partial orders $\dot Q_\al$, where $\dot Q_\al$ names 
$\P_\al$ if $\al$ is a regular uncountable cardinal, 
and $\dot Q_\al$ names the trivial forcing otherwise.
\end{defn}

Note that $\P_\ka$ will continue to denote the forcing to add a mangrove 
at $\ka$, whereas $R_\ka$ will now denote the iteration up to stage 
$\ka$ of such forcings.

\begin{thm}\label{GlobManPOWorks} 
If $V\sat\ZFC+\GCH$ and $G$ is $R$-generic over $V$, then
\begin{multline*}
V[G]\sat\ZFC+\GCH+\\
\text{``There is a mangrove at }\ka
\text{ for every regular cardinal }\ka\text{.''}
\end{multline*}
\end{thm}
\begin{proof}
For each regular cardinal $\ka$, we have that $\P_\ka$ is $\ka$-closed
and $\ka^+$-cc 
(Propositions \ref{MangroveForcingKClosed} and \ref{MorForcK+cc}),
and moreover $|\P_\ka|=(\ka^+)^\ka=\ka^+$.
As mentioned in Subsection~\ref{MangroveForcCardPres}, this implies
that forcing with $\P_\ka$ preserves cardinals, and a nice names argument
gives that the GCH is preserved as well.

Applying the Factor Lemma,
we may factorise $R$ as \mbox{$R_\ka*\dot R^\ka$} for any $\ka$, with
$R^\ka$ forced to be $\ka$-closed.  
As discussed in the proof of Lemma \ref{GCHForcTame}, 
this implies that $R$ is tame.

We prove that cardinals and the GCH are preserved by induction on the
iteration length.  Successor stages are immediate from the observations
about $\P_\ka$ above.

At singular limit cardinal $\la$ stages, 
we observe by factorising $R_\la$ as $R_{\ka^+}*R^{\ka^+}$ that
each cardinal $\ka<\la$ is preserved, whence $\la$ itself is also
preserved, and moreover the GCH is preserved below $\la$.  
The iteration $R_\la$ has a 
dense suborder of size $\prod_{\ka<\la}\ka^+=\la^+$, 
so cardinals greater than $\la^+$ are preserved, and the GCH holds 
at and above $(\la^+)^V$.
Finally, taking unions over $\ka<\la$ of nice names for 
subset of $\ka$, we obtain $\la^+$ nice names for subsets
of $\la$, so there is a surjection from $(\la^+)^V$ to $(2^\la)^{V[G_\la]}$.
Hence, $\la^+$ is also preserved, as is the GCH at $\ka$.

At regular limit stages $\la$, where we take direct limits, the argument is
similar but easier, since $R_\la$ will have a dense suborder of size
$\sum_{\ka<\la}\ka^+=\la$.  Indeed if $\la$ is Mahlo, we even have that
$R_\la$ satisfies the $\la$ chain condition --- see \cite{Bau:IF},
Corollary 2.4.  In any case, we have shown that every stage of the iteration
preserves cardinals and the GCH, and so since the closedness of $R^\ka$
increases with $\ka$, $R$ itself preserves cardinals and the GCH.

Showing that mangroves exist at every regular uncountable cardinal of 
the generic extension is now straightforward.
Given a regular uncountable $\ka$ and 
factorising $R$ as $R_\ka*\dot\P_\ka*\dot R^{\ka+1}$,
we have that there is a mangrove at $\ka$ after the first
$\ka+1$ stages of the iteration.
Now the statement ``$M$ is a mangrove at $\ka$'' is absolute
for models containing $M$, $\ka$ as a regular cardinal, and $\ka^+$ as
its successor.  Thus, since the final part $R^{\ka+1}$ of the
iteration is $\ka^+$-closed, the mangrove at $\ka$ will remain 
a mangrove after this final part of the iteration,
and we are done.
\end{proof}

\begin{thm}\label{GMornSuperstrongPres}
Suppose that $V\sat\ZFC+\GCH$, and $G$ is $R$-generic over $V$.
Then for any $n\in\omega+1$ and any cardinal $\ka$, if
$$
V\sat\ka\text{ is $n$-superstrong}
$$
then
$$
V[G]\sat\ka\text{ is $n$-superstrong}.
$$
Similarly, if $\ka$ is hyperstrong in $V$, then $\ka$ is hyperstrong
in $V[G]$.
\end{thm}
\begin{proof}

As mentioned above, the case of $1$-superstrongs is proven in
\cite{SDF:LCL}.  
The proofs in the other cases also follow exactly as in
that paper, with the simple modification that because of our
homogeneity theorem, a generic $G'\subset V[G]$ with $V[G']=V[G]$
may be chosen containing any given master condition.  
Now, we should be careful at this point, because for $n$-superstrength
with $n>1$, we will want to take a master condition at not just a
single iterand of the forcing, but many.  
Factorise $R$ as 
$$R_\ka*\dot\P_\ka*R^{(\ka,j(\ka))}*\cdots
*\P_{j^n(\ka)}*R^{j^n(\ka)+1},
$$
and let $G$ be $R$-generic over $V$. 
We wish to modify (within $V[G]$) the upper part of $G$ so as to contain
the master conditions for the lower part. 
There are no restrictions on $G_{j(\ka)}$.
We can modify $G$ so that $G(j(\ka))$ contains the master condition
from $G(\ka)$.
Now $R^{(j(\ka),j^2(\ka))}$ is $j(\ka)^+$-directed-closed
since each of its iterands is (see \cite{Bau:IF}, Theorem~2.7),
so we may take a single master condition 
in $R^{(j(\ka),j^2(\ka))}$ determined by $G^{(\ka,j(\ka)}$.
By arguments as in \cite{DoF:HIM1C},
our iteration $R^{(j(\ka),j^2(\ka))}$ is sufficiently homogeneous
to allow us to modify $G^{(j(\ka),j^2(\ka))}$, and hence $G$,
to lie below this 
single master condition.
Repeating all of this another $n-2$ times, 
we construct a $G'_{j^{n-1}(\ka)}\subset V[G]$ below all of the necessary 
master conditions to make the arguments of \cite{SDF:LCL} go through.
Of course, if $n=\omega$ the question arises as to whether this is legitimate,
but since an indirect limit is taken at stage $j^\omega(\ka)$
the methods of \cite{DoF:HIM1C}
are still applicable.

Therefore,
each $n$-superstrong or hyperstrong cardinal 
is preserved in \emph{every} generic extension by 
$R$, and so any generic extension $V[G]$ preserves all such cardinals.
\end{proof}

\section{1-extendible cardinals and morasses}\label{Mangrove1Extend}

In the case of forcing the GCH in Section \ref{GCH}, 
there was no major difficulty in proving that 1-extendible cardinals 
were preserved.  In particular, the fact that any embedding $j$ witnessing
1-extendibility lifted presented no great problems, 
essentially because the forcing affecting the domain of $j$
was small --- $P_\ka$ was a set in $H_{\ka^+}$.  To show that 
1-extendible cardinals are preserved while forcing morasses to exist,
we will have to work significantly harder, as the forcing of interest,
$R_{\ka+1}$,
will be class sized from the point of view of $H_{\ka^+}$.

\begin{thm}\label{1ExtMan}
Let $V$ be a model of $\ZFC$.
There is class generic extension $V[G]$ of $V$ satisfying the GCH,
having mangroves
at every regular cardinal in $V[G]$, and such that every
1-extendible cardinal of $V$ is 1-extendible in $V[G]$.
\end{thm}
\begin{proof}
We of course apply the two stage iteration $P*R$ of the GCH Partial Order from
Section \ref{GCH} followed by the Global Mangrove Partial Order above;
recall that whilst taking infinite iterations of class forcings is generally
problematic, there is no difficulty with 2-stage iterations.
We have already seen in Theorems \ref{GCH1Ext} and \ref{GlobManPOWorks} 
that if $G_g*G_m$ is $P*R$-generic over $V$, then $V[G_g*G_m]$ will
satisfy $\ZFC+\GCH$ and have mangroves at every regular cardinal.  
It therefore only remains to show that any 1-extendible cardinal $\ka$ of $V$
will remain 1-extendible in $V[G_g*G_m]$.

So suppose $\ka$ is 1-extendible in $V$, 
and let $j:H_{\ka^+}^{V}\to H_{\la^+}^V$ be an elementary embedding 
witnessing the 1-extendibility of $\ka$.
Let $j^*$ be the lift of $j$ to $H_{\ka^+}^{V[G_g]}$, 
as obtained in Theorem \ref{GCH1Ext}, so $j^*$ witnesses that
$\ka$ is 1-extendible in $V[G_g]$.
By Theorem \ref{GCH1ExtG} we may assume that $G_g$ is such that
$j^*$ is elementary from 
$\langle\HkaP^{V[G_g]},G_g(\ka)\rangle$ to 
$\langle\HlaP^{V[G_g]},G_g(\la)\rangle$ 
as $\LSTG$-structures, without changing $V[G_g]$.

Now the Global Mangrove Partial Order $R$ may be factorised as
$R\cong R_{\ka+1}*\dot R^{\ka+1}$, 
with $\dot R^{\ka+1}$ forced to be $\ka^+$-closed.
Hence, 
$$
H_{\ka^+}^{V[G_g*G_m]}=H_{\ka^+}^{V[G_g*G_{m,\ka+1}]},
$$ where
$G_{m,\ka+1}=G_m\restr R_{\ka+1}$.

To lift our embedding 
$$
j^*:\HkaP^{V[G_g]}\to\HlaP^{V[G_g]}
$$ 
to 
$$
j^{**}:\HkaP^{V[G_g*G_m]}\to\HlaP^{V[G_g*G_m]}, 
$$
we will first lift it to $\HkaP^{V[G_g*G_{m,\ka}]}$, 
and then again to $\HkaP^{V[G_g*G_{m,\ka+1}]}=\HkaP^{V[G_g*G_m]}$.

Lifting $j^*$ to 
${j^*}':\HkaP^{V[G_g*G_{m,\ka}]}\to\HlaP^{V[G_g*G_{m,\la}]}$
is like the case of forcing the GCH while preserving 1-extendible
cardinals (that is, Theorem~\ref{GCH1Ext}) and at first glance
seems to pose no problems.
In particular, 
$H_{\ka^+}^{V[G_g*G_{m,\ka}]}=H_{\ka^+}^{V[G_g]}[G_{m,\ka}]$,
and $j^*\restr R_\ka$ is the identity function, so the Lifting Lemma
gives us that 
${j^*}':\HkaP^{V[G_g]}[G_{m,\ka}]\to\HlaP^{V[G_g]}[G_{m,\la}]$
is elementary.
However, we actually want to demonstrate that ${j^*}'$ is 
elementary from 
$$\langle\HkaP^{V[G_g]}[G_{m,\ka}],G_g(\ka),G_{m,\ka}\rangle$$ to
$$\langle\HlaP^{V[G_g]}[G_{m,\la}],G_g(\la),G_{m,\la}\rangle.$$
The Lifting Lemma will certainly accommodate this if the forcing relation
is definable
for pre-existing class predicates $A$
(thus dealing with $G_g(\ka)$; $G_{m,\ka}$ is not a problem since
$R_\ka$ is still only set sized from the perspective of $\HkaP^{V[G_g]}$).
But we must be careful, as the definition of $\forces^*A$ as given
for example in \cite{SDF:FSCF}, page 34,
relies on our ground model having
a definable stratification into sets.

There are (at least) two ways to deal with this.
First, we could redefine $p\forces^*A(\si)$,
analogously to the definition of $p\forces^*G(\si)$ in the proof
of Lemma~\ref{GCHForcingHkaPDefinable}.  Specifically, it is not hard to
check that 
\begin{equation*}
p\forces^*A(\si)\quad\Iff\quad
\forall q\leq p\exists r\leq q\exists a\in A(r\forces^*\si=\check a)
\label{ClassForcADefn'}
\end{equation*}
is an appropriate definition.

Alternatively, we can give away the 
punch line\index{punch line} of
the proof that the next lift goes through, and observe that
our ground model $\langle\HkaP^{V[G_g]},G_g(\ka)\rangle$ 
\emph{does} have a definable stratification into sets.
Specifically, every element of $\HkaP^{V[G_g]}=\HkaP^{V[G_{g,\ka}]}$ 
is coded into the Cohen generic $G_g(\ka)$, and so
$\HkaP^{V[G_g]}=L_{\ka^+}[G_g(\ka)]$.
Hence, the $L_\al[G_g(\ka)]$ hierarchy is a stratification of $\HkaP^{V[G_g]}$
into sets, $p\forces A(\si)$ is definable, and the Lifting Lemma gives
that ${j^*}'$ is indeed elementary from 
$\langle\HkaP^{V[G_g]}[G_{m,\ka}],G_g(\ka),G_{m,\ka}\rangle$ to
$\langle\HlaP^{V[G_g]}[G_{m,\la}],G_g(\la),G_{m,\la}\rangle$.
That is, we have
$$
{j^*}':\langle\HkaP^{V[G_g*G_{m,\ka}]},G_g(\ka),G_{m,\ka}\rangle\to
\langle\HlaP^{V[G_g*G_{m,\la}]},G_g(\la),G_{m,\la}\rangle
$$
elementary.

We next wish to lift ${j^*}'$ to 
$j^{**}:
\HkaP^{V[G_g*G_{m,\ka+1}]}\to
\HlaP^{V[G_g*G_{m,\la+1}]}$.
Since $\P_\ka$ enjoys the $\ka^+$-cc, every $\P_\ka$-name in
$V[G_g*G_{m,\ka}]$
has an equivalent name in $\HkaP^{V[G_g*G_{m,\ka}]}$, so once again
this is the same as lifting to $\HkaP^{V[G_g*G_{m,\ka}]}[G_m(\ka)]$.
The forcing partial order $\P_\ka$ is class-sized from the point of view 
of $\HkaP^{V[G_g*G_{m,\ka}]}$, so we need a stratification of 
$\HkaP^{V[G_g*G_{m,\ka}]}$ into sets in $\HkaP^{V[G_g*G_{m,\ka}]}$
in order to deduce from Theorem~2.18 of \cite{SDF:FSCF} that the forcing
relation is definable, and 
subsequently apply the Lifting Lemma.
Note on the other hand that the pretameness of $\P_\ka$ is unproblematic 
for us: we know (from the $V$ perspective) that forcing with $\P_\ka$ 
over $\HkaP^{V[G_g*G_{m,\ka}]}$ will give us a model of $\ZFC^-$,
$\HkaP^{V[G_g*G_{m,\ka+1}]}$, so since pretameness is implied by
$\ZF^-$ preservation (Proposition~2.17 of \cite{SDF:FSCF}), 
$\P_\ka$ is pretame.  Alternatively, observe that 
the $\ka^+$-cc directly implies that $\P_\ka$ is pretame.

Now we have $\HkaP^{V[G_g]}=L_{\ka^+}[G_g(\ka)]$,
so 
$$
\HkaP^{V[G_g*G_{m,\ka}]}=\HkaP^{V[G_g]}[G_{m,\ka}]
=L_{\ka^+}[G_g(\ka),G_{m,\ka}].
$$
Thus, the $L_\al[G_g(\ka),G_{m,\ka}]$ hierarchy is a stratification of 
$\HkaP^{V[G_g*G_{m,\ka}]}$ into sets, which moreover is definable in 
$\langle\HkaP^{V[G_g*G_{m,\ka}]},G_g(\ka),G_{m,\ka}\rangle$.
So we can conclude that the forcing relation is
definable, and by the Lifting Lemma, if 
$${j^*}'``G_m(\ka)\subseteq G_m(\la),$$
then ${j^*}'$ lifts to an elementary 
$$
j^{**}:\HkaP^{V[G_g*G_{m,\ka}]}[G_m(\ka)]\to
\HlaP^{V[G_g*G_{m,\la}]}[G_m(\la)]
$$
as desired.
But now $|G_m(\ka)|=\ka^+$ and 
$G_m(\ka)$ is directed, so since $\P_\la$ is
$\la$-directed-closed (Lemma~\ref{MangroveForcingKClosed}),
there is a lower bound $p\in\P_\la$
for ${j^*}'``G_m(\ka)$.
By Theorem~\ref{MangroveForcingHomog}, there is a 
$G'_m(\la)\subset\HlaP^{V[G_g*G_{m,\la}]}[G_m(\la)]$ 
which is $\P_\la$-generic over $\HlaP^{V[G_g*G_{m,\la}]}$,
contains the master condition $p$, and such that
$$
\HlaP^{V[G_g*G_{m,\la}]}[G_m(\la)]= 
\HlaP^{V[G_g*G_{m,\la}]}[G'_m(\la)].
$$ 
Therefore, without changing the model $\HlaP^{V[G_g*G_{m,\la+1}]}$ we
may apply the Lifting Lemma, and obtain an elementary embedding
$$
j^{**}:\HkaP^{V[G_g*G_{m,\ka+1}]}\to\HlaP^{V[G_g*G_{m,\la+1}]}.
$$
As noted above, $\HkaP^{V[G_g*G_{m,\ka+1}]}=\HkaP^{V[G_g*G_m]}$, and
likewise with $\ka$ replaced by $\la$, so 
$j^{**}$ witnesses that $\ka$ is 1-extendible in $V[G_g*G_m]$.
\end{proof}

\section{Universal morasses}\label{univMorasses}

We now consider a notion extending that of a ($\ka, 1$)-morass, 
in which a predicate at the top of the morass, built up through the morass,
codes every subset of $\ka$.  
These are of interest to practitioners of pcf theory
for use in the construction of scales.


\begin{defn}\label{Aug+UMorass}
An \emph{augmented $\ka$-morass}\index{morass!augmented}
is a ($\ka, 1$)-morass $M$,
along with a set $A_x$ associated to each element $x$ of $\calS^M$, 
such that
\begin{enumerate}
\item for every $x\in\calS^M$, $A_x\subseteq o(x)$,
\item\label{AMorLevels} for $w,x\in\calS^M$ with $l(w)=l(x)$ and $o(w)<o(x)$, 
$A_w=A_x\cap o(w)$, and
\item\label{AMorPi} 
for $x,y\in\calS^M$ with $x\morl^My$, $A_x=(\pixy^M)^{-1}``A_y$.
\end{enumerate}
For any augmented morass $\langle M,A\rangle$, 
we denote by $A_\ka$\index{$A_\ka$} the set
$\bigcup_{\tau<\ka^+}A_{\langle\ka,\tau\rangle}\subseteq\ka$. 
A $\ka$-\emin{universal morass}\index{morass!$\ka$-universal} 
is an augmented $\ka$-morass
$\langle M,A\rangle$
such that $\Power(\ka)\subset L[A_\ka]$.
\end{defn}


It turns out that by modifying the mangrove forcing appropriately,
we may force a $\ka$-universal morass.  The idea is this: 
assuming we already have a subset $A'_\ka$ of $\ka^+$ such that 
$\Power(\ka)\subset L[A'_\ka]$, 
we may force an augmented morass 
with $A'_\ka$ at the top level.
However, there are new subsets of $\ka$ added by forcing the 
universal morass itself.  Thus, we also want the morass to be
self-coding, giving a subset of $\ka^+$ at the top level which 
encodes the morass itself.

For those familiar with Jensen's definition of a morass, 
as in \cite{Dev:Con} for example,
this will seem entirely straightforward:
at each node of the morass we can encode the morass below that node,
and then the morass maps will naturally preserve this information,
cohering at level $\ka$ of the morass to give a subset of $\ka^+$
encoding the morass itself.
However, in the definition of a morass that we have used
(Definition~\ref{morass}, taken from \cite{Vell:MDF}), 
there is no requirement that the morass maps are this ``nice''.
%
%
%
We are confident that 
by finding the right balance between Jensen's 
definition of a morass and Velleman's, this elegant solution might be made 
to go through.  
This will have the added benefit that attempts to make our suborder 
homogeneous will have one less problem to worry about: since
the augmentation of the morass for self-coding can be thought of
as occurring after the forcing takes place, it will not interfere with
any desired homogeneity.
On the other hand, new ideas will in any case be necessary to make our
suborder homogeneous while coding up all of the ground model subsets of $\ka$.

In the meantime, we use a different technique to
resolve the matter of self-coding.  
In essence, we shall force
our morass to be self-coding by treating the final code for the morass
much like the pre-existing subset $A'_\ka$.  
To do this, we place the restriction on our conditions that they are only
allowed to hit morass coding points at which membership or not
is decided by the condition itself.
The top level subset of the augmented morass that we force
will code the given $A'_\ka$ at 
successor ordinals,
and the morass itself will be coded at $0$ and limit ordinals,
so that there
will be no problems with $S^p$ being closed under ordinal successors and 
predecessors for conditions $p$.
For terminological convenience, we shall consider 0 to be a limit ordinal 
for the rest of this article.
Note that the function $\omega\cdot:\al\mapsto\omega\cdot\al$
is the increasing enumeration of all limit ordinals.

In defining our coding
we shall work a little harder than is necessary for our main goal,
in order to show how close one can get to homogeneity for our 
partial order, and where it can nevertheless fail.
For this, we wish to avoid coding redundant information.
It turns out that mangals will be helpful for this,
and moreover the following particular kind.

\begin{defn}\label{DetMan}
Given a morass or morass condition $M$,
a \emin{determined mangal} of $M$ is a mangal $\be$
of $M$ such that 
\begin{enumerate}
\item $\theta_\be$ 
is a limit ordinal, and
\item\label{DetManCov} for all $y\in M$ with $l(y)=\be$, 
there is a $\tau<\ka^+$ such that
$y\morl\langle\ka,\tau\rangle$. 
\end{enumerate}
\end{defn}

We will need the following basic result stemming from the morass axioms,
showing that no ``splitting'' occurs at limit stages of the tree.
\begin{lemma}\label{MorNoLimSplit}
Let $M$ be a morass or morass condition, and let $\be\leq\ka$ or $\la^M$
respectively.
Then 
for all $y\neq y'\in M$ with $l(y)=l(y')=\be$ and $o(y)+1<\theta_\be$ and
$o(y')+1<\theta_\be$,
there are
$x\neq x'\in M$ such that 
$x\closermorl y$, $x'\closermorl y'$, and $x$ and $x'$ are 
$\closermorl$-incomparable.
\end{lemma}
\begin{proof}
Take $y$ and $y'$ with $l(y)=l(y')=\be$ and $o(y')<o(y)$.
By M.5, there is an $x\morl y$ and a $\nu<o(x)$ such that
$\pixy(\nu)=o(y')$, whence from M.2 we have that 
$\langle l(x),\nu\rangle\morl y'$.  Taking 
$x'=\langle l(x),\nu\rangle$\mbox{$\morl$}$y'$, then, we are done.
\end{proof}

\begin{lemma}\label{DetManDetd}
Suppose $\be$ is a determined mangal of $M$.  Then the data
$\theta_\be$, the tree structure given by $\morl$ up to level $\be$ and 
directly from level $\be$ to
level $\ka$, and the corresponding maps $\pi$, are all determined by the
part of the morass below level $\be$ and the morass maps from such
lower nodes to level $\ka$.
\end{lemma}
\begin{proof}
First we reconstruct $\theta_\be$. 
Let $X\subseteq\ka^+$ be the set of $\tau\in\ka^+$ such that there is
an $x\in M$ with $l(x)<\be$ and $x\morl\langle\ka,\tau\rangle$, 
and define an equivalence relation $\sim$
on $X$ by
$$
\si\sim\tau\iff 
\forall x\big(l(x)<\be\implies
(x\morl\langle\ka,\si\rangle\iff
x\morl\langle\ka,\tau\rangle)\big).
$$
Then by Lemma~\ref{MorNoLimSplit}, 
$\theta_\be$ is in bijection with the set of $\sim$-equivalence classes,
and moreover thanks to M.2, $\theta_\be$ is 
the order type of the set of $\sim$-equivalence classes
ordered by least element.
For any $\ga<\theta_\be$ and $\tau\in\ka^+$, 
$\langle\be,\ga\rangle\morl\langle\ka,\tau\rangle$ if and only if
$\tau$ is in the $\ga$-th $\sim$-equivalence class of $X$, and similarly
for any $x\in M$ with $l(x)<\be$, $x\morl\langle\be,\ga\rangle$
if and only if there is a $\tau<\ka^+$ in the $\ga$-th equivalence class of
$X$ such that $x\closermorl\langle\ka,\tau\rangle$.

We also claim that the corresponding maps $\pi$ may be recovered.
For $x\morl\langle\be,\ga\rangle$ this is clear: 
if $\mu<o(x)$ and $\nu<\ga$,
then $\pi_{x\langle\be,\ga\rangle}(\mu)=\nu$ if and only if
for $\tau<\ka^+$ such that $x\morl\langle\ka,\tau\rangle$, 
$\pi_{x\langle\ka,\tau\rangle}(\mu)$ 
is in the $\nu$-th equivalence class of $X$.
The maps $\pi_{\langle\be,\ga\rangle\langle\ka,\tau\rangle}$ for
$\langle\be,\ga\rangle\morl\langle\ka,\tau\rangle$ may now be deduced
from the Commutativity morass axiom: for any $\nu<\ga$ we have by M.5
that there is an $x\morl\langle\be,\ga\rangle$ and a $\mu<o(x)$ such that
$\pi_{x\langle\be,\ga\rangle}(\mu)=\nu$, and then
$\pi_{\langle\be,\ga\rangle\langle\ka,\tau\rangle}(\nu)$ must be
$\pi_{x\langle\ka,\tau\rangle}(\mu)$.
\end{proof}

Because of this ``reconstructibility'' property, when we code up the
morass we need only flag the determined mangals as such, rather than 
all of the corresponding details.  Thus, when we iteratively construct
conditions that are sufficiently closed with respect to hitting the
morass coding points, we will be able to add on the necessary new top level
without adding condition-specific morass coding point information.

We are now ready to define our encoding for the morass.
The first $\ka\cdot\ka$ morass coding points will be used to encode
the values
$\tal$ for $\al<\ka$.
The next $\ka\cdot\ka$ will be used to encode the morass relation $\morl$
below level $\ka$,
and the next $\ka^4$ will be used to encode the maps
$\pixy$ below level $\ka$.  The rest of the morass coding points will be
used to encode the relation $x\morl y$ and maps $\pixy$ for $y$ with level
$\ka$.

\begin{defn}\label{MorCode}
Let $p$ be a morass condition, and for all $\ga\leq\la^p+1$ define
\begin{eqnarray*}
\ND(\ga)&=&\{\al<\ga\st\al\text{ is not a determined mangal of }p\}\\
\io(\ga)&=&\sum_{\al<\ga}\tal,\qquad\text{and}\\
\de(\ga)&=&\sum_{\al\in \ND(\ga)}\tal.
\end{eqnarray*}
Let $m:\io(\la^p+1)\to\calS^p$ $n:\de(\la^p+1)\to\calS^p$ 
be the enumeration in lexicographic order of the
nodes of $p$ with level in $\la^p+1$ and $\ND(\la^p+1)$
respectively.  
Let $b:\ka^+\times\ka^+\to\ka^+$ be the standard bijection
according to the max-lex order on $\ka^+\times\ka^+$.
Then the 
\emph{morass code for $p$}\index{morass code}
is a partial function 
$$c^p:\{\omega\cdot\al\st\al\in\ka^+\}\to 2$$
such that the following hold.
\begin{enumerate}
\item For all $\al\in\ND(\la^p+1)$ and $\zeta<\ka$,
$$
c^p(\omega\cdot(\ka\cdot\al+\ze))=\begin{cases}
1&\text{if }\ze<\tal\\
0&\text{otherwise.}
\end{cases}
$$
For every $\be$ a determined mangal of $p$ and $\ze<\ka$,
$c^p(\omega\cdot(\ka\cdot\be+\ze))=0$.
\item\label{MorCodemorl}
For $\al<\io(\la^p)$, let $\mu$ be the least determined mangal of $p$
greater than $l(m(\al))$ if such exists, and $\la^p+1$ otherwise.  
Then for all $\ze<\io(\mu)$,
$$
c^p(\omega\cdot(\ka^2+\ka\cdot\al+\zeta))=
\begin{cases}
1&\text{if }m(\al)\morl^p m(\zeta)\\
0&\text{otherwise.}
\end{cases}
$$
If $\mu$ is a determined mangal, then for all $\ze\geq\io(\mu)$,
$c^p(\omega\cdot(\ka^2+\ka\cdot\al+\zeta))=0$.

\item\label{MorCodepixy} For all $\al<\io(\la^p)$, 
all $\zeta<\io(\mu)$ with $\mu$ as in item \ref{MorCodemorl},
and all $\eta<\ka$,
\begin{multline*}
c^p(\omega\cdot(\ka^2\cdot 2+\ka^2\cdot\al+\ka\cdot\ze+\eta))\\
=
\begin{cases}
1&\text{if }m(\al)\morl^p m(\ze)\text{ and } 
\eta\in\pi_{m(\al)m(\ze)}``(o(m(\al))+1)\\
0&\text{otherwise.}
\end{cases}
\end{multline*}
If $\mu$ is a determined mangal, then for all $\ze\geq\io(\mu)$ and
$\eta<\ka$,
$$
c^p(\omega\cdot(\ka^2\cdot 2+\ka^2\cdot\al+\ka\cdot\ze+\eta))=0.
$$

\item For all $\tau\in S^p$, $\al<\de(\la^p+1)$, and $\nu<\ka^+$,
\begin{multline*}
c^p(\omega\cdot(\ka^3+\ka\cdot b((\tau,\nu))+\al))\\
=\begin{cases}
1&\text{if }n(\al)\morl^p\langle\ka,\tau\rangle\text{ and }
\nu\in\pi_{n(\al)\langle\ka,\tau\rangle}``(o(n(\al))+1)\\
0&\text{otherwise}
\end{cases}
\end{multline*}

\item At all other limit ordinals $c^p$ is undefined.
\end{enumerate}
\end{defn}

Clearly the entire condition $p$ may be reconstructed from the information
provided by $c^p$.
Also, our coding scheme is nice with respect to the mangrove $\leq$ relation:

\begin{lemma}
If $p$ and $q$ are morass conditions with $p\leq q$ in the mangrove forcing,
then $c^q\subseteq c^p$.
\end{lemma}
\begin{proof}
If new nodes, tree edges and attendant maps are added to extend $q$ to
$p$, these will be coded at limit ordinals on which $c^q$ was left undefined.
If new tree edges from nodes in $q$ to nodes at level $\ka$ are added, then
the node at level $\ka$ must be new (by injectivity of $f^p$) and again
the code is only extended at points where it was previously undefined.
\end{proof}

Now we must force in a way that respects $c^p$ and our predicate $A'_\ka$,
and moreover such that we will be able to respect $c^q$ for extensions $q$.
The modification to the quagmire forcing necessary to achieve this is given
by the following definition.  As usual, for any subset $X$ of $\ka^+$
we denote by $\chi_X$ the characteristic function of $X$.

\begin{defn}\label{UMorCond}
Let $A'$ be a subset of $\ka^+$ such that $\Power(\ka)\subset L[A']$,
and let $\chi_{A'}$ be its characteristic function.
Let $p$ be a morass condition.  Let 
$\chi^p_A:\ka^+\to\ka^+$ be the partial function defined by
$$
\chi_A^p(\tau)=\begin{cases}
\chi_{A'}(\nu)&\text{if }\tau=\nu+1\\
c^p(\tau)&\text{if }c^p\text{ is defined on }\tau
\end{cases}
$$
(and undefined otherwise), 
and let 
$A_\ka^p$ denote the subset of $\dom(\chi_A^p)$ 
with characteristic function $\chi_A^p$, 
The condition $p$ is a \emph{universal morass condition for $A'$}, or simply 
\emin{universal morass condition} when $A'$ is clear from the context,
if following properties hold.
\begin{enumerate}
\item\label{UMorCondSindomc} 
The set $S^p\cap\Lim$, where $\Lim$ denotes the class of limit ordinals,
is a subset of $\dom(c^p)$.
\item\label{UMorCondAxUnique} 
For all $x\in\calS^p$, if $\tau$ and $\nu$ in $S^p$ are such that
$x\morl^p\langle\ka,\tau\rangle$ and \mbox{$x\morl^p\langle\ka,\nu\rangle$,}
then 
$(\pi_{x\langle\ka,\tau\rangle}^p)^{-1}``(A_\ka^p\cap\tau)=
(\pi_{x\langle\ka,\nu\rangle}^p)^{-1}``(A_\ka^p\cap\nu)$.

\item\label{UMorCondAxDefined}
For all $x\in\calS^p$ with $l(x)\leq\la^p$, there is a $y\in p$ such that
$l(y)=\ka$ and $x\morl^py$.
\end{enumerate}
Let $U$ denote the set of all universal morass conditions for $A'$.  
The
\emph{Universal Morass Partial Order $\U_\ka$ for $A'$}
\index{Universal Morass Partial Order}
\index{partial order!Universal Morass}
is the suborder of $\P_\ka$ 
consisting of all universal morass conditions for $A'$,
that is, $\U_\ka=\langle U,\leq\rangle$ where $\leq$ is defined as for 
the mangrove forcing.
\end{defn}
As with $\P$, we will generally omit the subscript $\ka$ from $\U_\ka$.

We claim that forcing with $\U$ yields a universal morass.  
The corresponding set $A_\ka$ will be $\bigcup_{p\in G}A_\ka^p$,
and the other sets $A(x)$ for $x$ in the generic morass will
``filter down'' from this top level predicate:
if $x,y$ are nodes in the generic morass with $l(y)=\ka$, then 
$A(x)=\pixy^{-1}``(A_\ka\cap o(y))$.
Requirements \ref{UMorCondAxUnique} and \ref{UMorCondAxDefined}
ensure that $A(x)$ will be well-defined for all $x$ in the generic morass,
and the coherence properties necessary to have an augmented morass
will be automatic.

The verification that forcing with $\U$ still yields a morass 
is much as in Section~\ref{MangroveForcing},
with a little extra work to check that 
requirements \ref{UMorCondSindomc}--\ref{UMorCondAxDefined}
of Definition~\ref{UMorCond}
do not cause
problems.  We shall now walk through the modified versions
of the enunciations of Section~\ref{MangroveForcing}, 
only describing the necessary modifications to the old proofs.
In particular,
the fact that straightforward verifications go through will not 
be mentioned. \\


\noindent{\bf Lemma~\ref{NonTrivMorCond}$'$.}
\emph{There is a universal morass condition $q$ with $S^q=\omega$.}
\begin{proof}
The construction is exactly as for Lemma~\ref{NonTrivMorCond}.
Since $\langle0,0\rangle$ must be in $q$, $c^q(0)$ is defined 
(and equal to $1$).  
Every node of the $q$ constructed is below a unique node at level $\ka$,
and so $q$ satisfies all of the requirements 
of Definition~\ref{UMorCond}.
\end{proof}

In Subsection~\ref{muequiv} regarding $\mu$-equivalence, 
Lemmas~\ref{mcequivIncrStr} and 
\ref{MorCondCompatMangal} go through unchanged, 
but Lemma~\ref{MorCondChangeS} now fails in general.  
However, we do not require that lemma for our present purposes.

\noindent{\bf Proposition~\ref{MangroveForcingKClosed}$'$}
\emph{The poset $\U$ is $\ka$-directed-closed.}
\begin{proof}
As in the proof of Proposition~\ref{MangroveForcingKClosed}, given a
directed set $Y$ of conditions, we take the union of the conditions and
add a new top level if necessary.
The requirement of Definition~\ref{UMorCond}
hold below the new top level because they do for every member of $Y$,
and \ref{UMorCondAxUnique} and \ref{UMorCondAxDefined} hold at the 
new top level by its definition.
\end{proof}

We need to suitably strengthen the hypotheses of Lemma~\ref{ManForcDiffSCompat}
to make the proof go through in the context of $\U$.

\noindent{\bf Lemma~\ref{ManForcDiffSCompat}$'$}
\emph{Let $p$ and $q$ be universal morass conditions such that $\la^p=\la^q$
and $p\mcequiv{\la^p}q$. 
Suppose further that there is some $S^0\subset\ka^+$ such that both
$S^p$ and $S^q$ end-extend $S^0$, 
$\min(S^q\smallsetminus S^0)\geq\sup(S^p)$,
and $(f^p)^{-1}``A_\ka^p=(f^q)^{-1}``A_\ka^q$.  
Then $p$ and $q$ are compatible in $\U$.}
\begin{proof}
The added hypothesis ensures that requirement \ref{UMorCondAxUnique} holds
of the condition constructed
in the proof of Lemma~\ref{ManForcDiffSCompat},
and the other two requirements are immediate.
\end{proof}

\noindent{\bf Proposition~\ref{MorForcK+cc}$'$}
\emph{The poset $\U$ is $\ka^+$-cc.}
\begin{proof}
At the end of the $\Delta$-system argument for Proposition~\ref{MorForcK+cc},
we may add one extra stage of thinning, to ensure that 
$(f^p)^{-1}``A_\ka^p$ is the same for all of the conditions $p$ in the 
set in question, thus making Lemma~\ref{ManForcDiffSCompat}$'$ applicable.
\end{proof}

It will be convenient to consider the variants of 
Lemmas \ref{DsigmaDense} and \ref{MorCondVertExtn} in reverse order.

\noindent{\bf Lemma~\ref{MorCondVertExtn}$'$}
\emph{For any $p\neq\boldsymbol{1}\in\U$ there is a $q\leq p$ in $\U$ 
such that $S^q=S^p$,
$\theta_{\la^q}^q=\ot(S^q)$, and
$\la^q=\la^p+\omega^{\ot(S^q)}$.}
\begin{proof}
It is immediate that the $q$ constructed in the proof of 
Lemma~\ref{MorCondVertExtn} is a universal morass condition.
\end{proof}

\noindent{\bf Lemma~\ref{DsigmaDense}$'$}
\emph{For any $p\neq\boldsymbol{1}\in\U$ and any limit ordinal $\sigma<\ka^+$, there
is a $q\leq p$ in $\U$
such that $S^q\supseteq((\sigma+\omega)\smallsetminus\sigma)$.}
\begin{proof}
The only point of concern is to ensure that $\si\in\dom(c^q)$.
If $\si<\ka^3$, we may first
apply Lemma~\ref{MorCondVertExtn}$'$ repeatedly,
``gluing together'' with Proposition~\ref{MangroveForcingKClosed}$'$,
to obtain a $p'$ with sufficiently large $\la^{p'}$ and $\de(\la^{p'})$
that $c^{p'}$ is defined on $\si$.
The desired extension $q$ may then be 
constructed starting from $p'$ as in the proof of Lemma~\ref{DsigmaDense}.

If $\si=
\omega\cdot(\ka^3+\ka\cdot b((\tau,\nu))+\al)$
for some $\tau$, $\nu$ and $\al$, then $\tau\leq\si$.
We may again make a preliminary extension of $p$ to 
a $p'$ with $\al<\de(\la^{p'}+1)$. 
If $\tau=\si$, then the $q$ constructed in the proof of 
Lemma~\ref{DsigmaDense} satisfies $\si\in\dom(c^{q})$, and is thus
a universal morass condition.
Otherwise, by induction we may extend $p'$ to $p''$ with $\tau\in S^{p''}$,
which implies that $\si\in c^{p''}$, 
and so constructing $q$
from this starting point again makes it a
universal morass condition.
\end{proof}

We now have all the necessary ingredients for the following 
analogue of Theorem~\ref{limGMangrove},
which we believe merits its own number.

\begin{thm}\label{UMorForcWorks}
Forcing with $\U_\ka$ yields a $\ka$-universal morass.
\end{thm}
\begin{proof}
Let $G$ be $\U$-generic over $V$.
The morass part $M$ of the universal morass will simply be the direct limit of
the conditions in $G$, as was the case in our original mangrove forcing.
The proof that $M$ is indeed a morass proceeds exactly as for
Theorem~\ref{limGMangrove}.

As discussed earlier, we set $A_\ka=\bigcup_{p\in G}A_\ka^p$, and 
define $A(x)$ for every $x\in M$ by 
$$
A(x)=\begin{cases}
A_\ka\cap o(x)&\text{if }l(x)=\ka\\
\pixy^{-1}``A(y)&\text{if }x\morl y\text{ and }l(y)=\ka.
\end{cases}
$$
By requirements \ref{UMorCondAxUnique} and \ref{UMorCondAxDefined}
this uniquely defines $A(x)\subseteq o(x)$ for every $x\in M$.
Requirement \ref{AMorLevels} for an augmented morass is satisfied thanks
to morass axiom M.2 and our definition of $A(x)$, and
requirement \ref{AMorPi} for an augmented morass is immediate from
the Commutativity morass axiom.

Therefore, to prove that $\langle M,A\rangle$ so defined is a 
$\ka$-universal morass in $V[G]$, it only remains to show that 
every subset of $\ka$ in $V[G]$ lies in $L[A_\ka]$.
We have by assumption that $\Power(\ka)^V\subset L[A']$,
and $A'$ is coded into $A_\ka$ on successor ordinals.  
Now since $\U$ has the $\ka^+$ chain condition, 
every element of $\Power(\ka)^{V[G]}$ has a nice name in $H_{\ka^+}^V$, 
and therefore in $L[A_\ka]$.
But now the morass $M$ is entirely encoded into $A_\ka$ on limit ordinals,
and the generic $G$ may readily be reconstructed from the morass $M$
as the set of all universal morass conditions that are substructures of $M$.
Hence, $G\subset L[A_\ka]$, so within $L[A_\ka]$ 
we can correctly evaluate names, and therefore 
$\Power(\ka)^{V[G]}\subset L[A_\ka]$.  Overall then, we have that
$\langle M, A\rangle$ indeed constitutes a universal morass in $V[G]$.
\end{proof}

We may of course also take a reverse Easton iteration of these partial
orders $\U_\ka$ to force a $\ka$-universal morass to exist at every 
uncountable regular cardinal $\ka$.  For this, we will need to interleave
Cohen forcings to obtain the sets $A'$ with the forcing posets $\U_\ka$.

\begin{defn}\label{GlobUMorPO}
The \emin{Global Universal Morass Partial Order}
\index{partial order!Global Universal Morass} $W$ is the reverse
Easton iteration forcing with $\Fn(\ka^+,2,\ka^+)*\U_\ka$ at
uncountable regular cardinal stages $\ka$, and the trivial forcing 
elsewhere.  In the former case, 
the $A'$ from which the partial order $\U_\ka$ is defined is taken to be
the subset of $\ka^+$ whose characteristic function is the union of
the generic for the preceding $\Fn(\ka^+,2,\ka^+)$.
\end{defn}

\begin{thm}\label{GlobUMorPOWorks}
Forcing with the Global Universal Morass Partial Order
yields a universe in which the GCH holds and universal morasses
exist at every uncountable regular cardinal.
\end{thm}
\begin{proof}
This is a straightforward verification as for 
Theorem~\ref{GlobManPOWorks}.
Also note that we only need the GCH to hold below $\ka^+$ for the
lemmas leading up to Theorem~\ref{UMorForcWorks}, 
so we need not even start with a model of the GCH.
\end{proof}

Now we consider large cardinal preservation. 

\begin{thm}\label{GlobUMorPOnSuperstrong}
Forcing with the Global Universal Morass Partial Order $W$ preserves all
superstrong cardinals of the ground model.  Moreover,
given any $1$-extendible, hyperstrong or $n$-superstrong cardinal $\ka$ 
of the ground model
for $n\in\omega+1$,
the $W$-generic $G$ may be chosen so as to preserve the large cardinal
strength of $\ka$.  
\end{thm}
\begin{proof}
For $n$-superstrong and hyperstrong cardinals
this is exactly as for Theorem~12 of \cite{SDF:LCL}.
For 1-extendible cardinals, the situation is much as in Theorem \ref{1ExtMan}.
Given $j:H_{\ka^+}\to H_{\la^+}$, we first lift to 
$j':H_{\ka^+}^{V[G_\ka]}\to H_{\la^+}^{V[G_\la]}$ using the fact that
$W_\ka$ is a set forcing from the perspective of $H_{\ka^+}$.
Then we lift to 
$$
j'':\langle H_{\ka^+}^{V[G_\ka\ast G_{g,\ka}]},G_{g,\ka}\rangle\to
\langle H_{\la^+}^{V[G_\la\ast G_{g,\la}]},G_{g,\la}\rangle,
$$
where $G_{g,\ka}$ and $G_{g,\la}$ are the generics for $\Fn(\ka^+,2,\ka^+)$
and $\Fn(\la^+,2,\la^+)$ respectively, arguing as in the proof of 
Theorem \ref{GCH1ExtG}.  Finally, since $\U_\la$ is $\la$-directed-closed, 
we may choose the generic $G$ so that the corresponding part $G_{U,\la}$
lies below $j''``G_{U,\ka}$, and hence such that our embedding lifts to
$$
j^*:H_{\ka^+}^{V[G_\ka\ast G_{g,\ka}\ast G_{U,\ka}]}=H_{\ka^+}^{V[G]}\to
H_{\la^+}^{V[G_\la\ast G_{g,\ka}\ast G_{U,\la}]}=H_{\la^+}^{V[G]},
$$ 
as required.
\end{proof}

As a final point, 
we consider briefly the question of the homogeneity, or lack thereof,
of $\U$.  With our carefully defined coding scheme, we can find conditions
that seem to almost be sufficiently interchangeable to define an
automorphism of $\U$ as was done in Definition~\ref{Mangrovephipr}
for $\P$.
Let $\supp(c)$ denote the subset of $\dom(c)$ on which $c$ is not $0$.

\begin{prop}
Let $p$ and $q$ be universal morass conditions.
Then there are universal morass conditions
\mbox{$p'\leq p$} and $q'\leq q$ such that
$\la^{p'}=\la^{q'}$, $\theta_{\la^{p'}}^{p'}=\theta_{\la^{q'}}^{q'}$,
$\dom(c^{p'})=\dom(c^{q'}),$
$\la^{p'}$ is a determined mangal of both $p'$ and $q'$,
and
$\supp(c^{p'})\subseteq S^{p'}=S^{q'}\supseteq\supp(c^{q'}).$
\end{prop}
\begin{proof}
Having put in the effort to set up our coding carefully in
Definition~\ref{MorCode}, this is now relatively easy. 
Note that the fact that $\theta_{\la^{p'}}^{p'}=\theta_{\la^{q'}}^{q'}$
follows from $S^{p'}=S^{q'}$ and the requirement that every node
of a universal morass condition lie below a node at level $\ka$.
Now, we may construct sequences $p_i$ and $q_i$ for $i<\omega$ with
$p_0=p$, $q_0=q$, and for every $i\in\omega$,
$$
\supp(c^{p_i})\cup \supp(c^{q_i})\subseteq S^{p_{i+1}}\cap S^{q_{i+1}}
$$
and
$$
\dom(c^{p_i})\cup\dom(c^{q_i})\subseteq\dom(c^{p_{i+1}})\cap\dom(c^{q_{i+1}}).
$$
Note that this latter requirement implies that
$\la^{p_{i+1}}\geq\max(\la^{p_i},\la^{q_i})\leq\la^{q_{i+1}}$.
The construction is by combining Lemmas \ref{DsigmaDense}$'$ and
\ref{MorCondVertExtn}$'$ and Proposition \ref{MangroveForcingKClosed}$'$;
note that in the proofs of all three of these results, the condition
constructed has a determined mangrove as the top level, so in particular
parts \ref{MorCodemorl} and \ref{MorCodepixy} 
of the definition of the morass code
(Definition~\ref{MorCode}) do not pose a problem.
Applying Proposition~\ref{MangroveForcingKClosed}$'$ to obtain
lower bounds for the sequences $(p_i)_{i\in\omega}$
and $(q_i)_{i\in\omega}$ gives the desired $p'$ and $q'$ respectively,
since the new top level in each case is a determined mangal,
requiring no new 1s in the morass code.
\end{proof}

Unfortunately, this is still not sufficient to be able to define an
automorphism interchanging $p'$ and $q'$.  
The difficulty lies in the fact that
an extension $r$ of $p'$ may be able to ``see'' extra information about the
structure of $p'$ through the required agreement between
$(\pixy^r)^{-1}``A_\ka^r$ and
$(\pi_{xz}^r)^{-1}``A_\ka^r$ whenever $x\morl^ry$ and $x\morl^rz$ with
$l(y)=l(z)=\ka$.
Of course, this may contradict what ``happens in'' $q$, thus
blocking the existence of an automorphism of $\U$ interchanging $p'$ and $q'$.
So we do not have the requisite homogeneity needed to make Theorem
\ref{GlobUMorPOnSuperstrong} go through for a proper class of $n$-superstrong
or hyperstrong cardinals, and must leave this as an open question.
However, as mentioned above, it seems likely that a solution to this problem
will emerge with a suitable compromise
between the Jensen and Velleman definitions of a morass.


\bibliographystyle{asl}
\bibliography{lcg1m}

\end{document}